\documentclass[11pt]{amsart}
\usepackage{amsmath,amssymb,times,color}

\usepackage{fig4tex}
\definecolor{webred}{rgb}{0.75,0,0}
\definecolor{webgreen}{rgb}{0,0.75,0}
\usepackage[citecolor=webgreen,colorlinks=true,linkcolor=webred]{hyperref}
\usepackage{xspace}

\newtheorem{thm}{Theorem}[section]

\theoremstyle{definition}

\theoremstyle{remark}
\newtheorem{rem}[thm]{Remark}

\usepackage{stmaryrd} 
\usepackage{mathrsfs}

\oddsidemargin=\evensidemargin

\numberwithin{equation}{section}

\renewcommand{\Im}{\operatorname{Im}}

\newcommand{\fact}[1]{#1\mathpunct{}!}

\newcommand{\Div}{\operatorname{\mathrm{div}}}

\newcommand{\rot}{\operatorname{\mathrm{curl}}}

\DeclareMathOperator{\curl}{curl}
\DeclareMathOperator{\veccurl}{\mathbf{curl}}

\newcommand{\Ocal}{\mathcal{O}}
\newcommand{\Oc}{\mathcal{O}_{\mathrm{c}}}
\newcommand{\Oe}{\mathcal{O}_{\mathrm{e}}}

\newcommand  {\JJ}{\boldsymbol{\mathsf J}}
\newcommand{\Ec}{\EE_{\mathrm{c}}}
\newcommand{\Em}{\EE_{\mathrm{m}}}

\newcommand{\Ee}{\EE_{\mathrm{e}}}

 \newcommand{\xt}{{\boldsymbol{\mathsf{x}}_\mathsf{T}}}

\newcommand{\db}[1]{_{\raise-0.3ex\hbox{$\scriptstyle #1$}}}
\newcommand{\dd}[1]{_{\raise-1.5pt\hbox{$\scriptstyle #1$}}}

\newcommand{\di}{\displaystyle}

\newcommand{\dr}{{\rm d}}

\newcommand  {\C}{{\mathbb C}}

\newcommand  {\R}{{\mathbb R}}
\newcommand  {\T}{{\mathbb T}}

\newcommand {\Id}{\mathbb {I}}
\DeclareMathAlphabet{\mathbbold}{U}{bbold}{m}{n}

\renewcommand  {\H}{{\mathrm H}}
\newcommand{\Uc}{\UU_{\mathrm{c}}}
\newcommand{\Ue}{\UU_{\mathrm{e}}}

 \renewcommand{\T}{{\mathsf{T}}}

\newcommand{\BB}{\boldsymbol{\mathsf B}}

\newcommand  {\EE}{\boldsymbol{\mathsf E}}

\newcommand  {\HH}{\boldsymbol{\mathsf H}}

\newcommand  {\LL}{\boldsymbol{\mathsf L}}
\renewcommand  {\L}{{\mathrm L}}

\newcommand  {\RR}{\boldsymbol{\mathsf R}}

\newcommand  {\UU}{\boldsymbol{\mathsf U}}

\newcommand  {\VV}{\boldsymbol{\mathfrak H}}
\newcommand  {\WW}{{\boldsymbol{\mathfrak E}}}

\renewcommand  {\dd}{{\boldsymbol{\mathsf d}}}

\newcommand  {\eps}{\varepsilon}

\newcommand  {\nn}{\boldsymbol{\mathsf n}}

\newcommand  {\vv}{\boldsymbol{\mathsf v}}

\newcommand  {\xx}{\boldsymbol{\mathsf x}}

\newcommand{\Om}{\Omega}
\newcommand{\layer}{o}

\newcommand{\Otw}{\Omega_{\layer}^\eps}

\newcommand  {\cW}{{\mathfrak E}}

\newcommand  {\bH}{\mathbf{H}}

\newcommand  {\bL}{\mathbf{L}}

\newcommand  {\fke}{\mathfrak{e}}

\newcommand{\D}{\mathsf D}

\newcommand{\W}{\mathsf W}

\newcommand{\kc}{\kappa_{\mathrm{c}}}
\newcommand{\ke}{\kappa_{\mathrm{e}}}
\newcommand{\km}{\kappa_{\mathrm{m}}}

\newcommand{\muc}{\mu_{\mathrm{c}}}
\newcommand{\mue}{\mu_{\mathrm{e}}}
\newcommand{\mum}{\mu_{\mathrm{m}}}

\newcommand{\sigc}{\sigma_{\mathrm{c}}}
\newcommand{\epsc}{\epsilon_{\mathrm{c}}}
\newcommand{\sigm}{\sigma_{\mathrm{m}}}
\newcommand{\epsm}{\epsilon_{\mathrm{m}}}
\newcommand{\sige}{\sigma_{\mathrm{e}}}
\newcommand{\epse}{\epsilon_{\mathrm{e}}}

\newcommand{\eg}{\textit{e.g\mbox{.}}\xspace}

\newcommand{\Rd}{\color{red}}
\newcommand{\Bk}{\color{black}}

\begin{document}


\title[WKB Expansions and Asymptotic Models  for the Time-Harmonic Maxwell equations
]{
Wentzel--Kramers--Brillouin Expansions and  Generalized Impedance Transmission Conditions for 
 Thin-Layer Problems in Electromagnetism with Application to Biological Cells
}

\author[V. P\'eron]{Victor P\'eron}
 \address[V. P\'eron]{Universit\'e de Pau et des Pays de l'Adour, E2S UPPA, LMAP CNRS UMR 5142, 64013 Pau, France}
 \email[V. P\'eron]{victor.peron@univ-pau.fr}

\begin{abstract}
In this work we derive a WKB expansion for the electromagnetic fields solution of the time-harmonic Maxwell equations set in a domain with a thin layer. 
As a by-product of this expansion we obtain new second order asymptotic models with generalized impedance transmission conditions that turn out to depend on the mean curvature of the boundary of the subdomain surrounded by the thin layer. 
We show that these models can be easily integrated in finite element methods by developing mixed variational formulations.   One application of this work concerns the computation of the  electromagnetic field in biological cells.

\end{abstract}


\date{\today} 

\maketitle


\section{Introduction}

This paper presents new approximations of the electromagnetic fields for time-harmonic thin-layer transmission problems.   The simplified configuration is motivated for instance by the computation of the electromagnetic field in biological cells, see, {\it  e.g.},~\cite{DPP11,Durufle2014} where additional references may also be found.  This work is also motivated by other applications that involve biological cells and media with thin inclusions (see, {\it  e.g.},~\cite{ammari2014mathematical,cakoni2018nondestructive}) : biomedical applications,  nondestructive testing of delaminated interfaces between materials,  microwave imaging, radar applications,  geophysical applications$\cdots$.  

In this work we revisit an asymptotic method developed in recent works for solving scattering problems of  time-harmonic electromagnetic waves in domains with thin ayer~\cite{DPP11,Durufle2014}. In the problem of interest a main difficulty for the computation of the electromagnetic field with finite elements lies in the thinness of the layer since meshes with thin cells are required. To overcome this difficulty we develop a new asymptotic method based on a Wentzel--Kramers--Brillouin (WKB) expansion that makes possible to mimic the electromagnetic fields in the thin layer by approximate transmission conditions.  

The aim of this work is twofold. First, we derive a new WKB expansion   for the electric field in power series of a small parameter $\varepsilon$  which represents for instance the relative size of a cell membrane. We make explicit the first terms of this expansion and we infer new second order asymptotic models for the electric field with respect to $\varepsilon$. Second, we show how these new asymptotic models can be easily  integrated in finite element methods.  

There are several differences between this work and previous works~\cite{DPP11,Durufle2014} in which 
the authors derived $\varepsilon$-parameterizations for the electric field as $\EE^0 +\varepsilon
\EE^1$~\cite[Th. 2.9]{DPP11} and a second order model with Generalized Impedance Transmission Conditions (GITC)~\cite{Durufle2014} based on a multiscale expansion. However this multiscale expansion leads in general to an asymptotic modeling error  since the  considered  \emph{ansatz} for the electric field in the extra-cellular domain $\Oe^\eps$ (see Figure \ref{F1}) involves a sequence of asymptotics 
 providing a series that does not approach the electric field  inside the thin layer~\cite[Eq. (5.1a)-(5.1c)]{DPP11}.  
On the contrary, in this work tangential traces of the residues  are continuous across the boundaries of the thin layer by construction of the WKB expansion, see Section~\ref{S4.4}. That is why the WKB expansion provides in general sharper error estimates than the multiscale expansion derived in~\cite{DPP11,Durufle2014}. Therefore the new second order models provide also in general sharper error estimates than the second order GITC model derived in~\cite{Durufle2014}. 
Furthermore the new second order asymptotic models depend on the mean curvature of the interface $\Gamma$ between the subdomains $\Oc$ and $\Oe=\Ocal\setminus\overline{\Oc}$ (see Figure \ref{F1}). This is a main difference with the second order asymptotic models derived in the previous works~\cite{DPP11,Durufle2014}. This is also a difference  with many thin layer problems in electromagnetism (see, {\it  e.g.}, ~\cite{EN93,ammari1998effective,haddar2002stability,goffi2017transfer,goffi2020approximate,stupfel2021well} for the case of Maxwell's equations where additional references may also be found)  that lead in general to a second order impedance boundary condition for which the effect of the curvature does not appear.  We refer also the reader to the related but different work~\cite{khelifi2016small} where the authors derive an asymptotic expansion of the boundary perturbations of the electromagnetic fields resulting from small perturbations of the shape of a smooth inhomogeneity. \Bk


 There is a main advantage to use these new asymptotic models for numerical purposes with finite element computations since they do not require to mesh the thin layer and they can be simulated by using a mesh independent of  $\varepsilon$. On the contrary, there is a main difference between this work and the strategies developed for instance   by  Chun {\it  et al.}~\cite{CHH10} and Delourme {\it  et al.}~\cite{DHJ12} 
  to derive asymptotic models for thin-layer transmission problems in electromagnetism since  a different mesh is needed  for computing the numerical models for each value of $\varepsilon$,  {\it  e.g.}~\cite{Durufle2014}. \Bk  
 Furthermore,  as in~\cite{Durufle2014}, the new transmission conditions can also be easily integrated in finite element methods by   introducing an additional unknown. For this purpose we develop mixed variational formulations in Section~\ref{S3}.

This work is concerned essentially with the derivation of a WKB expansion  for the electric field and with the development of mixed variational formulations for asymptotic models. Theoretical aspects concerning the analysis of asymptotic models as well as  numerical aspects concerning  the performance of the new models in comparison with  previous works~\cite{CHH10,DHJ12,Durufle2014} will be devoted to further works. 

The outline of the paper proceeds as follow. In Section~\ref{S2}, we present the mathematical model for the electromagnetic fields. In Section~\ref{S2bis}, we present the first terms of the WKB expansion and  we compare the first terms of this WKB expansion  with the multi-scale expansion derived in~\cite{DPP11,Durufle2014}  (Section~\ref{Scompar}). Then we present a new second order GITC model (Section~\ref{Sgitc}) and an application to biological cells (Section~\ref{Se2.5}). In Section~\ref{S3}, we present  mixed variational formulations for GITC models. In Section~\ref{sec:Asymp}, we provide elements of proof for the WKB expansion.

\Bk

\section{The mathematical model
}
\label{S2}
After the introduction of notations in Section~\ref{Snot}, we introduce the mathematical model for the electric and magnetic fields, and the electric field formulation in  Section~\ref{MatMod}.

\newcommand{\Scal}{\mathcal{S}}
\newcommand{\Er}{\mathbb{R}}

\subsection{Notations}
\label{Snot}	
For any orientable smooth surface without boundary $\Scal$ of
$\Er^3$, the unit normal vector $\nn$ on $\Scal$ is outwardly
oriented from the interior domain enclosed by $\Scal$ towards the
outer domain. 

 We denote by $\veccurl_{\Scal}$ the tangential rotational operator (which applies to functions defined on $\Scal$) and $\rot_\Scal$ the surface rotational operator (which applies to vector fields)~\cite{N01} :
\begin{align*}
\forall \,f\in C^\infty(\Scal),\quad \veccurl_{\Scal} f&=\left(\nabla_\Scal f\right)\times  \nn\ ,
\\
\forall \,\vv\in \left(C^\infty(\Scal)\right)^3,\quad \rot_\Scal \vv&=
\Div_\Scal \left(\vv\times \nn\right) \ ,
\end{align*}
where $\nabla_\Scal$ and $\Div_\Scal$ are the tangential gradient and the 
surface divergence on $\Scal$,  respectively. 

We denote by $\bL_t^2(\Scal)$ the space of $\L^2$-integrable tangent vector fields on  $\Scal$ :
\begin{equation*}
\bL_t^2(\Scal) = \{\vv:\Scal  \rightarrow \R^{3} \, | \, 
\vv\in (\L^2(\Scal))^3 , \,  \vv\cdot \nn =0 \mbox{ on } \Scal \} \ ,
\end{equation*}
and, for any $s\in [-1,1]$, the Sobolev space of tangential vector fields of order $s$ 
on the surface $\Scal$ is defined as
\begin{equation*}
\T\H^{s}(\Scal) = \{\vv:\Scal  \rightarrow \R^{3}  \, | \, 
\vv\in (\H^{s}(\Scal))^3 , \,  \vv\cdot \nn =0 \mbox{ on } \Scal \} \ ,
\end{equation*}
where $\H^{s}(\Scal)$ is the classical Sobolev space of order  $s$ defined on  $\Scal$.  
In this framework we denote  by $\T\H(\Div_{\Scal},\Scal)$, 
$\T\H^{-1/2}(\Div_{\Scal},\Scal)$ and $\T\H(\curl_{\Scal},\Scal)$
the following spaces of tangent vector fields of the above operators  
$\Div_{\Scal}$  and  $\rot_\Scal$:
\begin{align*}
 &\T\H(\Div_{\Scal},\Scal)= \{\vv\in  \bL_t^2(\Scal), \, \Div_{\Scal} \vv\in  \L^2(\Scal)  \} \ ,
 \\
 &\T\H^{-1/2}(\Div_{\Scal},\Scal)= \{\vv\in \T\H^{-1/2}(\Scal), \, \Div_{\Scal} \vv\in  \H^{-1/2}(\Scal)  \} \ ,
 \\
&\T\H(\curl_{\Scal},\Scal)= \{\vv\in \bL_t^2(\Scal), \, \curl_{\Scal} \vv\in \L^2(\Scal)  \}\ . 
\end{align*}

Finally we denote by $\T\H^{-1/2}(\Div_{\Scal},\Scal,0) $ and $\T\H(\curl_{\Scal},\Scal,0)$ the following spaces of tangent vector fields 
\begin{align*}
  &\T\H^{-1/2}(\Div_{\Scal},\Scal,0)= \{\vv\in \T\H^{-1/2}(\Div_{\Scal},\Scal), \, \Div_{\Scal} \vv = 0 \,  \mbox{ on } \Scal \} \ ,
 \\
&\T\H(\curl_{\Scal},\Scal,0)= \{\vv\in \T\H(\curl_{\Scal},\Scal), \, \Div_{\Scal} \vv = 0 \, \mbox{ on } \Scal \}\ .
\end{align*}

For any  vector field $\vv$ defined in a neighborhood of $\Scal$ we denote by $\vv_{\T}|_{\Scal}$ the tangent
component of  $\vv$ :
$$\vv_{\T}|_{\Scal}= \nn\times \left( \vv|_{\Scal}\times\nn\right)\ .
$$
Finally we denote  by 
 $\left[\vv\right]_{\Scal}$ the jump of $\vv$ across $\Scal$ :
 $$
 \left[\vv\right]_{\Scal} = \vv|_{\Scal^+} - \vv|_{\Scal^-}\ ,
 $$
and we denote by  $\langle \vv|_{\Scal}\rangle$ the mean of $\vv$ across $\Scal$ :
$$
\langle \vv|_{\Scal}\rangle=\frac12 (\vv|_{\Scal^+}+ \vv|_{\Scal^-}) \ . 
$$

\subsection{Time-harmonic Maxwell equations in single cell}
\label{MatMod}
Biological cells consist of a cytoplasm surrounded by a thin layer, see, {\it  e.g.},~\cite{DPP11,Durufle2014}.  We denote by $\Ocal$ the three-dimensional domain of interest which is composed of the outer cell medium and the cell. Let us denote by $\Oc$ the cell cytoplasm, and by $\Ocal^\eps_{\mathsf{m}}$ the cell
membrane surrounding $\Oc$, whose thickness is constant and
denoted by $\eps$. Assuming, without loss of generality, that the
domain $\Oc$ is independent of $\eps$, the extracellular domain is
then $\eps$--dependent. We denote it by $\Oe^\eps$, in a such way that (see Figure \ref{F1})
:
$$\Ocal =\Oc \cup\overline{\Ocal^\eps_{\mathsf{m}}}\cup\Ocal^\eps_{\mathsf{e}}.$$
The boundary of the cytoplasm is the smooth surface denoted by
$\Gamma$ while $\Gamma^\eps$ is the cell boundary, {\it i.e.}
$\Gamma^\eps$ is  the
boundary of $\overline{\Oc}\cup\Ocal^\eps_{\mathsf{m}}$.

\begin{figure}[ht]
\begin{center}
\psset(width=2)
\figinit{1.3pt}
\figpt 1:(-100,20)
\figpt 2:(-30,70)\figpt 3:(50,50)
\figpt 4:(80,0)\figpt 5:(30,-40)
\figpt 6:(-30,0)\figpt 7:(-80,-20)
\figpt 8:(-51,-12)\figpt 9:(-20,50)\figpt 10:(50,20)
\figpt 12:(10,10)\figpt 13:(0,-10)
\figpt 14:(-20,15)%
\figpt 22:(-92,15.5)\figpt 23:(-27.2,60.8)
\figpt 24:(46.3,41.5)\figpt 25:(70.8,2.8)
\figpt 26:(29,-31)\figpt 27:(-27.2,9.2)\figpt 28:(-77,-10.8)
\figpt 36:(-130,40)\figpt 37:(-40,100)
\figpt 38:(80,70)\figpt 39:(120,0)
\figpt 40:(60,-60)\figpt 41:(-20,-70)\figpt 42:(-110,-50)
\figpt 16:(-70,20)\figpt 15:(0,30)
\figpt 17:(40,10) \figpt 18:(-60,-20)
\figpt 21:(50,65) \figpt 31:(30,-25)  \figpt 34:(-32,65) 
\figpt 35:(-80,-15.4) \figpt 43:(-30,-40) \figpt 44:(-90,80)
\figpt 19:(61,65) \figpt 20:(3,27)  \figpt 30:(42.6,33)  \figpt 32:(61.6,5.6) \figpt 33:(61.6,6.1)
\figpt 45:(20,20) \figpt 46:(20,15) \figpt 47:(10,63) \figpt 48:(20,59) \figpt 49:(-20,90) 
\figpt 50:(-20,80) 
\psbeginfig{}
\psset arrowhead(fillmode=no,length=3)\psarrow[2,23]
\psset arrowhead(fillmode=no,length=3)\psarrow[23,2]
\pscurve[1,2,3,4,5,6,7,1,2,3]
\pssetfillmode{no}\pssetgray{0}
\pscurve[22,23,24,25,26,27,28,22,23,24]
\pscurve[36,37,38,39,40,41,42,36,37,38]
\pssetfillmode{yes}\pssetgray{5.0}
\psendfig
\figvisu{\figBoxA}{}
{
\figwritew 15: $\Oc$(18pt)
\figwritec [34]{$\varepsilon$}
\figwritec [8]{$\Gamma^{\varepsilon}$} 
\figwritec [31]{$\Gamma$}
\figwritec [35]{$\Ocal^\eps_{\mathsf{m}}$}
\figwritec [43]{$\Oe^\varepsilon$}
\figwritec [46]{($\muc$, $\sigc$, $\epsc$)}
\figwritec [49]{($\mue$, $\sige$, $\epse$)}
\figsetmark{$\figBullet$}
}
\centerline{\box\figBoxA}
 \caption{A cross-section of the domain $\Ocal$ and its subdomains $\Oc$, $\Ocal^\eps_{\mathsf{m}}$, $\Oe^\eps$}
\label{F1}
\end{center}
\end{figure}

The electromagnetic properties of $\Ocal$ are given by the following
piecewise-constant functions $\underline{\mu}$,  $\underline{\epsilon}$, and $\underline{\sigma}$  corresponding to the magnetic
permeability, the electrical permittivity, and the conductivity of $\Ocal$, respectively :
  $$
  \underline{\mu}=
  \begin{cases}
    \muc,\quad \text{in $\Oc$},\\
    \mum,\quad \text{in $\Ocal^\eps_{\mathsf{m}}$},\\
    \mue,\quad \text{in $\Oe^\eps$},
  \end{cases}\quad
  \underline{\epsilon}=
  \begin{cases}
    \epsc,\quad \text{in $\Oc$},\\
    \epsm,\quad \text{in $\Ocal^\eps_{\mathsf{m}}$},\\
    \epse,\quad \text{in $\Oe^\eps$},
  \end{cases}\quad 
  \underline{\sigma}=
  \begin{cases}
    \sigc,\quad \text{in $\Oc$},\\
    \sigm,\quad \text{in $\Ocal^\eps_{\mathsf{m}}$},\\
    \sige,\quad \text{in $\Oe^\eps$},
  \end{cases}.  $$
Let us denote by $\JJ$ the time-harmonic  current source and let
$\omega$ be the frequency. 
For the sake of simplicity, we assume that $\JJ$ is  smooth, $\Div\JJ\in \L^{2}(\Ocal)$, and that $\JJ$ is supported in $\Oe^\eps$ and it vanishes in a neighborhood of the cell membrane. 
Maxwell's equations link the electric field $\EE$ and the magnetic
field $\HH$, through  Faraday's and Amp\`ere's laws in $\Ocal$ : 
$$
\curl\EE^\eps - i \omega \underline{\mu} \HH^\eps = 0 \quad \mbox{and} \quad \curl\HH^\eps + \left(i \omega \underline{\epsilon}-\underline{\sigma} \right) \EE^\eps = \JJ\quad \mbox{in} \quad 
\Ocal\ .
$$
 We complement this problem with a  Silver-M\" uller boundary condition set on
 $\partial\Ocal$.
In what follows we denote by $\underline\kappa$ the complex wave number given by
$$ \forall \xx\in\Ocal,\quad \underline{\kappa}^2 (\xx)= \omega^2 \underline{\mu}(\xx)
\left(\underline{\epsilon}(\xx)+i \frac{\underline{\sigma}(\xx)}{\omega}
\right),\quad \Im(\underline{\kappa}(\xx))\geqslant 0 \ .
$$ 
Then Maxwell's system of first order partial differential equations can be
 reduced to the following second-order differential equation 
\begin{subequations}
  \label{eq:Eeeps}
  \begin{align}
&\curl\curl \EE^\eps - \underline \kappa^2  \EE^\eps = i\omega\underline{\mu} \JJ\quad \mbox{in} \quad \Oc\cup\Ocal^\eps_{\mathsf{m}}\cup\Oe^\eps\ ,
\intertext{
with the following transmission conditions across $\Gamma$ and
$\Gamma^\eps$
}   
& \Ee^\eps \times \nn|_{\Gamma^\eps} \; = \; \Em^\eps \times \nn|_{\Gamma^\eps}, \quad
 \frac{1}{\mue} \curl{\Ee^\eps} \times \nn|_{\Gamma^\eps} \; = \; \frac{1}{\mum} \curl{\Em^\eps} \times \nn|_{\Gamma^\eps} ,\\
 &\Ec^\eps \times \nn|_{\Gamma} \; = \; \Em^\eps \times \nn|_{\Gamma}, \quad
\frac{1}{\muc} \curl{\Ec^\eps} \times \nn|_{\Gamma} \; = \;
\frac{1}{\mum} \curl{\Em^\eps} \times \nn|_{\Gamma},
\intertext{  
where $\Ee^\eps$, $\Em^\eps$, and $\Ec^\eps$ denote  the 
restrictions of $\EE^\eps$ to the subdomains $\Oe^\eps$,
$\Ocal^\eps_{\mathsf{m}}$, and $\Oc$,  respectively. Finally, the boundary condition is given as 
}&  \curl{\EE^\eps} \times \nn - i \ke\, \nn \times \EE^\eps \times \nn = 0 \quad\mbox{on} \quad \partial\Ocal . 
  \end{align}
\end{subequations}

\section{First terms of a WKB expansion and asymptotic models with generalized impedance transmission conditions
}
\label{S2bis}

In this section,  we present  the first terms of a WKB expansion (Section~\ref{midfreq}) and asymptotic models with generalized impedance transmission conditions (Section~\ref{Sgitc}), with application to biological cells in  Section~\ref{Se2.5}.

\subsection{First terms of a WKB expansion}
\label{midfreq}

In this work we derive a WKB expansion for the electric field $\EE^\eps$ \eqref{eq:Eeeps}  as follow (compare with~\cite[Eq. (5.1)]{Durufle2014}) 
\begin{subequations}
\label{eq:ME}
\begin{align}
\Ec^\eps(\xx)&
\approx \di\sum_{j\geqslant 0} \eps^{j} \Ec^{j} (\xx)
,  \quad \text{in}\quad \Oc  \ ,  \\
\label{eq:ME2}
\Em^\eps(\xx)&
\approx  \di\sum_{j\geqslant 0} \eps^{j} \left(\Ee^{j}(\xx)+\Em^j\left(\xt,\frac{x_3}{\eps}\right) \right)
, 
\quad \text{in}\quad  \Ocal^\eps_{\mathsf{m}}  \ ,  \\ 
\Ee^\eps(\xx)&
\approx  \di\sum_{j\geqslant 0} \eps^{j} \Ee^{j} (\xx)
,  \quad \text{in}\quad \Oe^\eps
 \ , 
\end{align}
\end{subequations}
where the terms $\Ec^{j}$, $\Ee^{j}$ and $\Em^j$  are independent of $\eps$ and defined in $\eps$-independent domains $\Oc$, $\Oe=\Ocal\setminus\overline\Oc$ and  $\Gamma\times [0,1]$, respectively ; in \eqref{eq:ME2} $\xt$ denotes tangential coordinates on $\Gamma$ and  the variable $x_3 \in(0,\eps)$ is the Euclidean distance to $\Gamma$. 

For such a purpose it is convenient to introduce the electromagnetic
properties of the ``background''  problem, {\it i.e} the domain $\Ocal$  without   the
membrane:
   $$
  {\mu}=
  \begin{cases}
    \muc,\quad \text{in $\Oc$},\\
    \mue,\quad \text{in $\Oe$},
  \end{cases}\quad
  {\epsilon}=
  \begin{cases}
    \epsc,\quad \text{in $\Oc$},\\
        \epse,\quad \text{in $\Oe$},
  \end{cases}\quad 
  {\sigma}=
  \begin{cases}
    \sigc,\quad \text{in $\Oc$},\\
    \sige,\quad \text{in $\Oe$},
  \end{cases}, $$
and we define similarly $\kappa$ as
$${\kappa}=
  \begin{cases}
    \kc,\quad \text{in $\Oc$},\\
    \ke,\quad \text{in $\Oe$} ,
  \end{cases} .
  $$

\subsubsection*{Terms of order $0$ of the WKB expansion  \eqref{eq:ME}}
The term  $\EE^0=(\Ec^{0}, \Ee^{0})$ of order $0$  defined  in the subdomains  $\Oc$ and $\Oe$  solves the following transmission problem : 
\begin{subequations}\label{0E0}
\begin{align}
 &\curl\curl \Ec^{0}-\kc^2 \Ec^{0}=0\ ,\,\text{in   $\Oc$},\\	
 &\curl\curl \Ee^{0}-\ke^2 \Ee^{0}=i \omega \mue \JJ\ ,\,\text{in  $\Oe$},
 \intertext{with the transmission conditions:}
 & \Ee^{0}\times \nn|_{\Gamma^+}= \Ec^{0}\times \nn |_{\Gamma^-}, 
 \quad
  \frac{1}{\mue} \curl{\Ee^0} \times \nn|_{\Gamma^+} \; = \; \frac{1}{\muc} \curl\Ec^0 \times \nn|_{\Gamma^-},
 \\
  & \left(i \omega \epse-\sige \right) \Ee^{0}\cdot \nn |_{\Gamma^+}= \left(i \omega \epsc-\sigc \right) \Ec^{0}\cdot \nn |_{\Gamma^-}, \quad 
 \intertext{and the Silver-M\" uller condition}
 & \curl{\Ee^0} \times \nn - i \ke\, \nn \times \Ee^0 \times \nn = 0 \quad\mbox{on} \quad \partial\Ocal . 
\end{align} 
\end{subequations} 
The term  $\Em^0$  of order $0$  defined in  $\Gamma\times [0,1]$ depends only on the variable $\xt\in\Gamma$ and we have  : 
\begin{equation}
\label{EEm0}
\Em^0\left(\xt\right)= 
\left(\dfrac{i \omega \epse-\sige}{i \omega \epsm-\sigm}-1\right)\Ee^{0}\cdot \nn |_{\Gamma^+}\, \nn \ .
\end{equation}

\subsubsection*{Terms of order $1$ of the WKB expansion  \eqref{eq:ME}}
The term  $\EE^1=(\Ec^{1}, \Ee^{1})$   of order $1$  defined in the subdomains  $\Oc$ and $\Oe$   solves the following problem
\begin{subequations}
\label{0E1} 
\begin{align}
&\curl\curl \Ec^{1}-\kappa^2 \Ec^{1}=0\
,\,\text{in $\Oc\cup\Oe$}\ ,
\\
&\curl{\Ee^1} \times \nn - i \ke\, \nn \times \Ee^1\times \nn = 0 \quad\mbox{on} \quad \partial\Ocal ,
 \intertext{
with the following transmission conditions on $\Gamma$ (compare with~\cite[Eq. (3.4c)-(3.4d)]{Durufle2014}) : 
}
 \begin{split}
&\Ee^{1}\times \nn|_{\Gamma^+}= \Ec^{1}\times \nn |_{\Gamma^-}
+ \left(\frac{\mum}{\km^2}-\frac{\mue}{\ke^2}\right)\veccurl_{\Gamma} \curl_{\Gamma} (\frac1{\muc}\curl \Ec^{0})_{\T}
\\&-\left(\mum-\mue\right)(\frac1{\muc}\curl \Ec^{0})_{\T} 
\ ,
\label{TC1}
\end{split}
\intertext{ and } 
\begin{split}
&\frac{1}{\mue}\left(\curl \Ee^{1}\times \nn\right)|_{\Gamma^+}=\frac{1}{\muc}\left(\curl \Ec^{1}\times \nn
\right)|_{\Gamma^-}
\\& +\left(\frac{1}{\mum}-\frac{1}{\mue}\right)
\left(\ke^2 \  (\Ee^{0})_{\T}
+
\veccurl_\Gamma\rot_\Gamma   (\Ee^{0})_{\T}
\right)
\\  &+ \frac{2}{\mum} \frac{\mue}{\ke^2} \left(1-\dfrac{i \omega\epse-\sige}{i \omega  \epsm-\sigm}\right) 
\nn\times\veccurl_{\Gamma} \left( \mathcal{H}\curl_{\Gamma} (\frac1{\mue}\curl \Ee^{0})_{\T}\right)
\ .
\label{TC2}
\end{split}
\end{align}
\end{subequations}
In \eqref{TC2} $\mathcal{H}$ is the {\it mean curvature} of the surface $\Gamma$, $\veccurl_{\Gamma}$ is the tangential rotational operator and $\rot_\Gamma$ is the surface rotational operator, see Section \ref{Snot}. 

The term  $\Em^1$  of order $1$  defined in  $\Gamma\times [0,1]$  is defined as :
\begin{multline}
\label{EEm1}
\Em^1\left(\xt,\frac{x_3}{\eps}\right)=
 \left(a +  b \frac{x_3}{\eps} \right) \,  \nn 
\\ +  \left(\frac{x_3}{\eps}-1\right)  
\left( \nn\times\left(\frac{\mum}{\km^2}-\frac{\mue}{\ke^2}\right)\veccurl_{\Gamma} \curl_{\Gamma} (\frac1{\muc}\curl \Ec^{0})_{\T}
+\left( \frac{\mum}{\mue} - 1 \right)\left(\curl \Ee^{0}\times \nn\right) |_{\Gamma^+} 
\right)
\ ,
\end{multline}
 for almost any $(\xt,x_3)\in\Gamma\times (0,\eps)$, where functions  $a$ and $b$ are defined on $\Gamma$ by 
\begin{equation}
\label{E1aE1b}
a=\left(\dfrac{i \omega\epsc-\sigc}{i \omega  \epsm-\sigm}-1\right)\Ec^{1} \cdot \nn |_{\Gamma^{-}} - \Ee^{1}\cdot \nn |_{\Gamma^{+}}\ , 
 \quad b= 2\mathcal{H} \left(\dfrac{i \omega\epse-\sige}{i \omega  \epsm-\sigm}-1\right)\Ee^{0}\cdot \nn |_{\Gamma^{+}}  \ .
\end{equation}

\subsubsection{Comparison with the multi-scale expansion derived in~\cite{DPP11,Durufle2014}}
\label{Scompar}

The first term $\EE^0$ \eqref{0E0}  coincides with the ``background'' model derived in~\cite{DPP11}. Note also that the transmission condition \eqref{TC1} for the next term $\EE^1$ coincides with the transmission condition derived in~\cite[Eq. (3.4c)]{Durufle2014}. However there is a main difference between the new transmission condition \eqref{TC2} and the transmission condition derived in~\cite[Eq. (3.4d)]{Durufle2014} which writes (by abusing the notation we denote also by $\EE^{1}$  below  the second term of the multi-scale expansion derived in~\cite{DPP11,Durufle2014})
\begin{multline*}
\frac{1}{\mue}\left(\curl \EE^{1}\times \nn\right)|_{\Gamma^+}=
\frac{1}{\muc}\left(\curl \EE^{1}\times \nn
\right)|_{\Gamma^-}
\\
-\left(\frac{\km^2}{\mum}-
\frac{\ke^2}{\mue}\right) \EE_{\T}^{0}|_{\Gamma^+}  
+\left(\frac{1}{\mum}-
\frac{1}{\mue}\right)
\veccurl_\Gamma\rot_\Gamma
\EE_{\T}^{0}|_{\Gamma^+} \ .
\end{multline*}
Note in particular that in difference with the previous work~\cite[Eq. (3.4)]{Durufle2014} the new term $\EE^1$ defined by \eqref{0E1} depends on the mean curvature $ \mathcal{H}$ of the interface $\Gamma$. Furthermore all the terms describing the thin layer in \eqref{eq:ME2} are different from the corresponding asymptotics of the multiscale expansion derived in~\cite{DPP11,Durufle2014}. 

\subsection{GITC model}  
\label{Sgitc}

The electric field $\EE_{[1]}^\eps$ which approximates $\EE^\eps$ at the second order  with respect to the parameter $\eps$ is obtained by solving the problem : 
\begin{subequations}
\label{GITC}
  \begin{align}
&\curl\curl \EE_{[1]}^\eps-\kappa_{\mathrm{e}}^2 \EE_{[1]}^\eps=i \omega \mue\JJ\
,\,\text{in  $\Oe$ },\\
&\curl\curl \EE_{[1]}^\eps-\kappa_{\mathrm{c}}^2 \EE_{[1]}^\eps=0 \
,\,\text{in $\Oc$},\\
 &\curl{\EE_{[1]}^\eps} \times \nn - i \ke\, \nn \times \EE_{[1]}^\eps\times \nn = 0 \quad\mbox{on} \quad \partial\Ocal ,
\intertext{with the following transmission conditions, called 
generalized impedance transmission conditions (GITC) of order 2 (compare with~\cite[Eq. (4.2d)-(4.2e)]{Durufle2014}):}
 &\left[ \nn \times \EE_{[1]}^\eps \right]_{\Gamma} =\varepsilon \left(- A \, \veccurl_{\Gamma} \rot_{\Gamma} \left< \frac{1}{\mu}(\curl\EE_{[1]}^\eps|_{\Gamma} )_\T\right> \;+ B \left< \frac{1}{\mu}(\curl\EE_{[1]}^\eps|_{\Gamma} )_\T\right> \right) \ ,
 \label{GITC_4}
  \medskip \\
\begin{split}
 \label{GITC_5}
 & \left[ \frac{\nn}{{\mu}} \times \curl{\EE_{[1]}^\eps} \right]_{\Gamma}  =  \varepsilon \left(\,- C \, \veccurl_{\Gamma} \rot_{\Gamma} \left< (\EE_{[1]}^\eps)_{\T} \right> \; + \; D \left< (\EE_{[1]}^\eps)_{\T} \right>\right) 
\\
&  \qquad\qquad \qquad- \eps E  \; \nn\times\veccurl_{\Gamma} \mathcal{H} \curl_{\Gamma}  \left< \frac{1}{\mu}(\curl\EE_{[1]}^\eps|_{\Gamma} )_\T  \right> \ ,
 \end{split}
\end{align}
\end{subequations}
where constants $A,B, C, D, E$ are defined in \eqref{ABCD} 
\begin{multline}
\label{ABCD}
A \; = \;  \frac{\mum}{\km^2} - \frac{\mue}{\ke^2}  \medskip, 
\quad B\; = \; \mum -\mue \medskip ,
\quad C\; = \; \frac{1}{\mum} - \frac{1}{ \mue} \medskip , 
\quad  D \; = \;\left( \frac{1}{\mue} - \frac{1}{\mum}\right)\ke^2  \medskip , 
\\ 
\quad  E \; =  \frac{2}{\mum}\frac{\mue}{\ke^2} \left(1-\dfrac{i \omega\epse-\sige}{i \omega  \epsm-\sigm}\right) 
\	 .
\end{multline}

\begin{rem} 
In the framework of Section~\ref{VarFramwk} a proof of 
uniform estimates for $\EE_{[1]}^\eps\in\mathbf{V}$, where the functional space  $\mathbf{V}$ is defined in Section~\ref{S3} by \eqref{EV}, is non trivial since there is a lack of control of the divergence of the electric fields in  $\mathbf{V}$ that prevents to obtain a compact embedding of this space in $\bL^{2}(\Ocal)$. 
One way to prove such results is to exhibit a Helmholtz decomposition of  $\mathbf{V}$ (see, {\it e.g.},~\cite{DHJ12}). However one can not adapt straightforwardly the proof of~\cite[Prop. 9]{DHJ12} to obtain a Helmholtz decomposition of  $\mathbf{V}$  since in difference with the framework of Delourme {\it  et al.}~\cite{DHJ12} as well as with the previous work~\cite[Eq. (4.2)]{Durufle2014}, there is an extra-term  $\nn\times\veccurl_{\Gamma} \mathcal{H} \curl_{\Gamma}  \left< \frac{1}{\mu}(\curl\EE_{[1]}^\eps|_{\Gamma} )_\T  \right>$ which appears in the new  transmission condition \eqref{GITC_5}. 
\end{rem}

\subsection{Application to biological cells}
\label{Se2.5}

 Biological cells are characterized by constant magnetic permeabilities $\mum=\mue$. This allows to simplify the above asymptotic models.

\subsubsection{Two first orders of the asymptotic expansion}
 
In the framework above the WKB expansion for the electric field $\EE^\eps$ 
still writes by abusing the notations : 
\begin{align*}
\EE^\eps&
\approx{\EE}^{0}+\eps{\EE}^{1} +\cdots ,  \quad
 \text{in}\quad \Oc\cup\Oe^\eps \\
\EE^\eps&
\approx 
{\EE}^{0} +\Em^0\left(\xt \right) + \eps \left({\EE}^{1}+\Em^1\left(\xt,\frac{x_3}{\eps}\right) \right) +\cdots, 
\quad \text{in}\quad  \Ocal^\eps_{\mathsf{m}}
 \ ,
\end{align*}
where the first term $\EE^0$ still solves Problem \eqref{0E0} and $\Em^0\left(\xt \right)$ is defined by \eqref{EEm0}, whereas the term $\EE^1$ solves  the following simplified transmission problem in comparison with  Problem \eqref{0E1}  
\begin{subequations}
\begin{align}
&\curl\curl \EE^{1}-\kappa^2 \EE^{1}=0\
,\,\text{in $\Oc\cup\Oe$},\\
 & \curl{\EE^1} \times \nn - i \ke\, \nn \times \EE^1\times \nn = 0 \quad\mbox{on} \quad \partial\Ocal ,
 \intertext{with the following transmission conditions on $\Gamma$} 
&\EE^{1}\times \nn|_{\Gamma^+}= \EE^{1}\times \nn |_{\Gamma^-} +\left(\frac{\mum}{\km^2}-\frac{\mue}{\ke^2}\right)\veccurl_{\Gamma} \curl_{\Gamma} (\frac1{\muc}\curl \Ec^{0})_{\T}
\ ,
\\
\begin{split}
&\frac{1}{\mue}\left(\curl \EE^{1}\times \nn\right)|_{\Gamma^+}=\frac{1}{\muc}\left(\curl \EE^{1}\times \nn
\right)|_{\Gamma^-} 
\\
&\qquad \qquad\qquad+ \frac{2}{\ke^2} \left(1-\dfrac{i \omega\epse-\sige}{i \omega  \epsm-\sigm}\right) 
\nn\times\veccurl_{\Gamma} \left( \mathcal{H}\curl_{\Gamma} (\frac1{\mue}\curl \Ee^{0})_{\T}\right)
\ .
\end{split}
\end{align}
\end{subequations}
Furthermore, the term $\Em^1\left(\xt,\frac{x_3}{\eps}\right)$  defined by   \eqref{EEm1} is also simplified as follows 
\begin{equation*}
\Em^1\left(\xt,\frac{x_3}{\eps}\right)= \left(\frac{x_3}{\eps}-1\right)  
\left( \nn\times\left(\frac{\mum}{\km^2}-\frac{\mue}{\ke^2}\right)\veccurl_{\Gamma} \curl_{\Gamma} (\frac1{\muc}\curl \Ec^{0})_{\T}\right) 
+\left(a +  b \frac{x_3}{\eps} \right) \,  \nn \ ,
\end{equation*}
 where functions  $a$ and $b$ are  defined on $\Gamma$ by \eqref{E1aE1b}.

\subsubsection{GITC model} 
\label{S2.5.2}

Constants $B,C$ and $D$ in~\eqref{ABCD}  are equal to zero. Hence, in the framework above applied to biological cells, the GITC model $\EE_{[1]}^\eps$~\eqref{GITC} is reduced to :
\begin{subequations}\label{GITC_cell}
  \begin{align}
&\curl\curl \EE_{[1]}^\eps-\kappa_{\mathrm{e}}^2 \EE_{[1]}^\eps=i \omega \mue\JJ\
,\,\text{in  $\Oe$ },
\label{GITC_cell1}
\\
&\curl\curl \EE_{[1]}^\eps-\kappa_{\mathrm{c}}^2 \EE_{[1]}^\eps=0 \
,\,\text{in $\Oc$},\\
 &\curl{\EE_{[1]}^\eps} \times \nn - i \ke\, \nn \times \EE_{[1]}^\eps\times \nn = 0 \quad\mbox{on} \quad \partial\Ocal ,
\intertext{with the following GITC  of order 2 :}
 &
 \left[ \nn \times \EE_{[1]}^\eps \right]_{\Gamma} = - \varepsilon A  \, \veccurl_{\Gamma} \rot_{\Gamma}\left< \frac{1}{\mu}(\curl\EE_{[1]}^\eps|_{\Gamma} )_\T\right>  \ ,
 \label{GITC_cell4}
 \medskip \\
  \label{GITC_cell5}
 & \left[ \frac{\nn}{{\mu}} \times \curl{\EE_{[1]}^\eps} \right]_{\Gamma}  = -  \varepsilon E  \; \nn\times\veccurl_{\Gamma} \mathcal{H} \curl_{\Gamma}  \left< \frac{1}{\mu}(\curl\EE_{[1]}^\eps|_{\Gamma} )_\T\right> 
 \ .
\end{align}
\end{subequations}

\section{Mixed variational formulations for GITC models}
\label{S3}

Since the transmission conditions~\eqref{GITC_4}-\eqref{GITC_5} and \eqref{GITC_cell4}-\eqref{GITC_cell5} are of Wentzel types,  we introduce an additional unknown $\lambda$ defined as
$$
\lambda=\left< \frac{1}{\mu}(\curl\EE_{[1]}^\eps|_{\Gamma} )_\T\right>\ , 
$$
in order to incorporate the transmission conditions in mixed variational formulations~\cite{Braess:2007,Durufle2014}. In this framework, conditions~\eqref{GITC_4}-\eqref{GITC_5} write as  
\begin{subequations}\label{GITCbis}
  \begin{align}
 &\left[ \nn \times \EE_{[1]}^\eps \right]_{\Gamma} =\varepsilon \left(- A \, \veccurl_{\Gamma} \rot_{\Gamma} \lambda \;+ B \lambda \right) \ ,
  \medskip \\
 \begin{split}
 & \left[ \frac{\nn}{{\mu}} \times \curl{\EE_{[1]}^\eps} \right]_{\Gamma}  =  \varepsilon \left(\,- C \, \veccurl_{\Gamma} \rot_{\Gamma} \left< (\EE_{[1]}^\eps)_{\T} \right> \; + \; D \left< (\EE_{[1]}^\eps)_{\T} \right> \right) 
\\
&  \qquad\qquad \qquad- \eps E  \; \nn\times\veccurl_{\Gamma} \mathcal{H} \curl_{\Gamma} \lambda	 \ ,
 \end{split}
\end{align}
\end{subequations}
and conditions~\eqref{GITC_cell4}-\eqref{GITC_cell5} write as  
\begin{subequations}\label{GITCter}
  \begin{align}
\begin{split}
 & \left[ \nn \times \EE_{[1]}^\eps \right]_{\Gamma} =  - \varepsilon A  \, \veccurl_{\Gamma} \rot_{\Gamma} \lambda \ ,   
 \medskip \\
 & \left[ \frac{\nn}{{\mu}} \times \curl{\EE_{[1]}^\eps} \right]_{\Gamma}  = -  \varepsilon E  \; \nn\times\veccurl_{\Gamma} \mathcal{H} \curl_{\Gamma}  \lambda 
 \ .
  \end{split}
\end{align}
\end{subequations}

\subsection{Variational frameworks}  
\label{VarFramwk}

In this section we introduce variational frameworks for the GITC models~\eqref{GITC} and \eqref{GITC_cell}. 

\subsubsection*{Framework for the GITC model~\eqref{GITC}}
First we introduce a variational framework for the GITC model~\eqref{GITC} where the transmission conditions are written with the additional unknown $\lambda$ as in \eqref{GITCbis}. For this purpose note that $\T\H(\curl_\Gamma,\Gamma)\cap \T\H(\Div_\Gamma,\Gamma)$ coincides with the Sobolev space $\T\H^1(\Gamma)$.  Then the functional spaces associated with $\EE_{[1]}^\eps$ and $\lambda$ are denoted by $\mathbf{V}$ and $ \W$, respectively, and are defined as 
\begin{subequations}
\label{EV}
\begin{align}
\begin{split}
& \mathbf{V}= 
\Bigl\{\EE \in \bL^2(\Ocal), \, \curl \Ec\in \bL^2(\Oc), \, 
\curl \Ee\in \bL^2(\Oe), \,   
\\
& 
\qquad \left< \EE_{\T}|_\Gamma \right> \in \T\H^1(\Gamma) , \,
\EE\times \nn\in \bL_t^2(\partial \Ocal) \Bigr\}\ ,
\end{split}
\\
&
\W=\T\H(\curl_{\Gamma},\Gamma).
\end{align}
\end{subequations}

\vspace{0.2cm}

\newcommand{\Gcal}{\mathcal{G}}
\subsubsection*{Framework for the GITC model~\eqref{GITC_cell}}
Second we introduce a variational framework for the GITC model~\eqref{GITC_cell}  where the transmission conditions are written with the additional unknown $\lambda$ as in \eqref{GITCter}. For this purpose we define the functional space $\W_{0}$ by 
\begin{equation}
\label{EW0}
\W_{0}=\T\H(\curl_{\Gamma},\Gamma,0) \ , 
\end{equation}
and let  $\mathcal{G}$ be the operator defined from  $\T\H^{-1/2}(\Div_{\Gamma},\Gamma,0) $ onto $\W_{0}$
 by 
\begin{align*}
   \text{for any 
 $g\in\T\H^{-1/2}(\Div_{\Gamma},\Gamma,0)$},\quad
&\Gcal(g)=\lambda ,\quad\text{
where $\lambda$ satisfies }
\notag\\
& A \, \vec{\curl}_{\Gamma} \rot_{\Gamma} \lambda 
= g
\quad \mbox{on} \quad \Gamma \ .
\end{align*}  
Then the operator $\mathcal{G}$ is invertible from  $\T\H^{-1/2}(\Div_{\Gamma},\Gamma,0) $ onto $\W_{0}$,
see~\cite{Durufle2014}. Finally the functional spaces associated with $\EE_{[1]}^\eps$ and $\lambda$ for the GITC model~\eqref{GITC_cell} written with the additional unknown $\lambda$ are $\mathbf{V}_{0}$ and $ \W_{0}$, respectively, where  $\mathbf{V}_{0}$ is defined as 
\begin{equation}
\label{EV0}
\begin{split}
& \mathbf{V}_{0}= 
\Bigl\{\EE \in \bL^2(\Ocal), \, \curl \Ec\in \bL^2(\Oc), \, 
\curl \Ee\in \bL^2(\Oe), \,   
\\
& 
\left< \EE_{\T}|_\Gamma \right> \in  \T\H(\Div_\Gamma,\Gamma), \,
\EE\times \nn\in \bL_t^2(\partial \Ocal) \Bigr\}\ \ ,
\end{split}
\end{equation}
and  $\W_{0}$ is defined as \eqref{EW0}. 

We are in position now to write mixed variational formulations  for both GITC models~\eqref{GITC} and \eqref{GITC_cell}.

\subsection{A mixed variational formulation
}  
\label{S4}

Assume that $\mum\neq\mue$ (i.e. $B\neq0$). 
The mixed variational formulation for the GITC model~\eqref{GITC} writes  (compare with~\cite[Eq. (4.8)]{Durufle2014}) : 

Find $(\EE_{[1]}^\eps, \lambda) \in  \mathbf{V}\times \W$ such that for any $(\UU,\xi) \in  \mathbf{V}\times \W$, 
\\[4pt] 
\begin{subequations}
\label{varformGlobal}
\begin{align}\label{varformGlobalE}
 \begin{split}
& 
\int_{\Oc} \frac{1}{\muc} \curl \EE_{[1]}^\eps \cdot \curl \overline{\Uc}\, \dr \xx 
+\int_{\Oe} \frac{1}{\mue} \curl \EE_{[1]}^\eps \cdot \curl \overline{\Ue}\, \dr \xx 
- \int_{\mathcal{O}} \frac{\kappa^2}{\mu} \EE_{[1]}^\eps \cdot \overline{\UU}\, \dr \xx 
\\
&- i \frac{\ke }{\mue}  \int_{\partial\mathcal{O}} \EE_{[1]}^\eps \times \nn \cdot \overline{\UU} \times \nn \, \dr s
 - \int_{\Gamma}  \lambda \cdot \left[ \, \overline{\UU} \times \nn\, \right]_{\Gamma} \medskip \, \dr s
 +  \varepsilon \int_{\Gamma} \, E\, \mathcal{H} \rot_{\Gamma}\lambda \rot_{\Gamma}  \left< \overline{\UU\times \nn} \right>\, \dr s
 \\
&+ \varepsilon \int_{\Gamma} \, C\, \rot_{\Gamma}  \left< (\EE_{[1]}^\eps)_{\T} \right> \rot_{\Gamma}  \left< \overline{\UU_{\T}} \right>\, \dr s
- \varepsilon \int_{\Gamma}D \left< (\EE_{[1]}^\eps)_{\T} \right>\cdot \left< \overline{\UU_{\T}} \right>\, \dr s
\\
&\; = \; i \omega \int_{\Oe} \JJ \cdot \overline{\Ue}\, \dr \xx \ ,
\end{split} 
\intertext{and }
\label{varformGlobalLambda}
& \int_{\Gamma} \left[ \nn \times \EE_{[1]}^\eps \right]_{\Gamma} \cdot \overline{\xi} \, \dr s+ \varepsilon \int_{\Gamma} A \rot_{\Gamma}  \lambda
 \rot_{\Gamma}  \overline{\xi}  \, \dr s - \varepsilon \int_{\Gamma} B\lambda \cdot \overline{\xi} \, \dr s \; = \; 0\ .
\end{align}
\end{subequations} \Bk
\begin{rem}
By abusing the notation 
we identify in \eqref{varformGlobalE}  the  integral $\int_{\Gamma}  \lambda \cdot \left[ \, \overline{\UU} \times \nn\, \right]_{\Gamma} \medskip \, \dr s$ with the duality product between the Hilbert spaces $\T\H(\curl_{\Gamma},\Gamma)$ and $\T\H^{-1/2}(\Div_{\Gamma},\Gamma)$. 
\end{rem}

\subsection{Application to biological cells}  
\label{S4.2}

In the framework of biological cells (Section~\ref{Se2.5}), the mixed variational formulation for the GITC model~\eqref{GITC_cell} writes : 

Find $(\EE_{[1]}^\eps, \lambda) \in  \mathbf{V}_{0}\times \W_{0}$ such that for any $(\UU,\xi) \in  \mathbf{V}_{0}\times \W_{0}$, 
\\[4pt] 
\begin{subequations}
\label{varformGlobal-cell}
\begin{align}\label{varformGlobalE-cell}
 \begin{split}
& 
\int_{\Oc} \frac{1}{\muc} \curl \EE_{[1]}^\eps \cdot \curl \overline{\Uc}\, \dr \xx 
+\int_{\Oe} \frac{1}{\mue} \curl \EE_{[1]}^\eps \cdot \curl \overline{\Ue}\, \dr \xx 
- \int_{\mathcal{O}} \frac{\kappa^2}{\mu} \EE_{[1]}^\eps \cdot \overline{\UU}\, \dr \xx 
\\
&- i \frac{\ke }{\mue}  \int_{\partial O} \EE_{[1]}^\eps \times \nn \cdot \overline{\UU} \times \nn \, \dr s
 - \int_{\Gamma}  \lambda \cdot \left[ \, \overline{\UU} \times \nn\, \right]_{\Gamma} \medskip \, \dr s
+  \varepsilon \int_{\Gamma} \, E\, \mathcal{H} \rot_{\Gamma}\lambda \rot_{\Gamma}  \left< \overline{\UU\times \nn} \right>\, \dr s  
\\
& \; = \; i \omega \int_{\Oe} \JJ \cdot \overline{\Ue}\, \dr \xx \ ,
\end{split}
\intertext{and }
\label{varformGlobalLambda-cell}
& \int_{\Gamma} \left[ \nn \times \EE_{[1]}^\eps \right]_{\Gamma} \cdot \overline{\xi} \, \dr s+ \varepsilon \int_{\Gamma} A \rot_{\Gamma}  \lambda
 \rot_{\Gamma}  \overline{\xi}  \, \dr s 
 = \; 0\ .
\end{align}
\end{subequations}
\Bk
We remind that constants $A, E$ are given by \eqref{ABCD} with $\mum=\mue$. 

\section{Derivation of a WKB expansion 
 for the electric field}
\label{sec:Asymp}

 In this section we denote by $\Omega$, $\Omega_{-}$, $\Omega^\eps_{\layer}$, and $\Omega^\eps_{+}$  the domain and  subdomains of interest instead of $\Ocal$, $\Oc$, $\Ocal^\eps_{\mathsf{m}}$, and $\Oe^\eps$,  respectively, and we denote the material properties by :
  $$
  \underline{\mu}=
  \begin{cases}
    \mu_{-},\quad \text{in $\Omega_{-}$},\\
    \mu_{\layer},\quad \text{in $\Omega^\eps_{\layer}$},\\
    \mu_{+},\quad \text{in $\Omega^\eps_{+}$},
  \end{cases}\quad
  \underline{\epsilon}=
  \begin{cases}
    \epsilon_{-},\quad \text{in $\Omega_{-}$},\\
     \epsilon_{\layer},\quad \text{in $\Omega^\eps_{\layer}$},\\
     \epsilon_{+},\quad \text{in $\Omega^\eps_{+}$},
  \end{cases}\quad 
  \underline{\sigma}=
  \begin{cases}
    \sigma_{-},\quad \text{in $\Omega_{-}$},\\
    \sigma_{\layer},\quad \text{in $\Omega^\eps_{\layer}$},\\
     \sigma_{+},\quad \text{in $\Omega^\eps_{+}$},
  \end{cases}.  $$

 The derivation of the generalized impedance transmission conditions in Sect.~\ref{Sgitc} is based on a WKB expansion 
 for the electric field~$\EE^\eps$~\eqref{eq:Eeeps} inside and outside the thin layer. More precisely, we search $\EE^\eps$ as the asymptotic expansion

 \begin{subequations}
\label{eq:MEbis}
\begin{align}
\EE_{-}^\eps (\xx) & \approx \di\sum_{j \geqslant 0} \eps^{j} \EE^{-}_{j}(\xx) 
,  \quad \text{in}\quad \Omega_{-}  \ ,  
\\[1.0ex]
\EE_{\layer}^\eps (\xx)  & \approx  \di\sum_{j \geqslant 0}\eps^{j} \left(\EE_{j}^{+}(\xx)+\WW_j \left(y_{\alpha},\frac{h}{\eps}\right) \right)
, 
\quad \text{in}\quad   \Omega^\eps_{\layer} \ ,  
\\[1.0ex] 
\EE_{+}^\eps  (\xx) &
 \approx  \di\sum_{j \geqslant 0} \eps^{j} \EE_{j}^{+}  (\xx) 
,  \quad \text{in}\quad  \Omega_{+}^\eps
 \ .
\end{align}
\end{subequations}
Here, $\left(y_{\alpha},h\right)$ is a 
{\em normal coordinate system} (see, \eg,~\cite[App. A.1]{CDFP11}) to the surface $\Gamma$ in the thin layer $\Otw$
where $y_{\alpha}$, $\alpha=1,2$ are tangential coordinates on $\Gamma$ and $h\in (0,\eps)$ is the normal coordinate to $\Gamma$.

In this section, we derive the first terms of this expansion step by step 
as well as their governing equations, having in mind that
the electric field~$\EE^\eps$ satisfies the following Maxwell transmission problem \Bk
\begin{subequations}
\label{eq:EME}
\begin{align}
  \label{eq:EME:1-}
   \rot\rot  \EE_-^{\eps} -\; \kappa_{-}^2\; \EE_-^{\eps}  &= 0  &&\mbox{in}\quad  \Omega_{-}\ ,
\\[1.0ex]
  \label{eq:EME:1+}
   \rot\rot  \EE_+^{\eps} -\; \kappa_+^2\; \EE_+^{\eps}  &= i \omega\mu_{+}\JJ   &&\mbox{in}\quad  \Omega^\eps_{+}\ ,
\\[1.0ex]
  \label{eq:EME:2}
  \rot\rot  \EE_\layer^{\eps} - (\kappa^\eps_{\layer})^2 \EE_\layer^{\eps} &= 0   &&\mbox{in}\quad \Omega^\eps_{\layer}\ ,
\\[0.5ex]
  \label{eq:EME:3}
 \dfrac1{\mu_{\pm}} \rot \EE_\pm^{\eps}\times\nn &=  \dfrac1{\mu_{\layer}}\rot \EE_\layer^{\eps}\times\nn &&\mbox{on} \quad\Gamma\cup\Gamma^{\eps}\ ,
\\[0.5ex]
  \label{eq:EME:4}
 \EE_\pm^{\eps}\times\nn &= \EE_\layer^{\eps}\times\nn  &&\mbox{on}\quad \Gamma\cup\Gamma^{\eps}\ ,
\\[1.0ex]
  \label{eq:EME:5}
\curl{\EE_{+}^{\eps}} \times \nn - i \kappa_{+}\, \nn \times \EE_{+}^{\eps} \times \nn &= 0 && \mbox{on} \quad \partial\Om\ .
\end{align}
\end{subequations}

Furthermore, according to Amp\`ere's law and by assumption on $\JJ$ we have : 
\begin{equation}
\label{eq:EME:div}
\Div\left(i \omega \underline{\epsilon}-\underline{\sigma} \right) \EE^\eps = \Div\JJ  \qquad \mbox{in} \quad \L^{2}(\Omega) \ . 
\end{equation}
We deduce the extra transmission conditions  
\begin{equation}
\label{eq:EMEbis}
 \left(i \omega \epsilon_{\pm}- \sigma_{\pm} \right)   \EE_\pm^{\eps}\cdot\nn =  \left(i \omega \epsilon_{\layer}
 -\sigma_{\layer} \right)  \EE_\layer^{\eps}\cdot\nn \qquad \mbox{on} \quad\Gamma\cup\Gamma^{\eps}\ .
\end{equation}



In Section~\ref{AH1}, we expand the "electric" Maxwell operators inside the thin layer $\Omega^\eps_{\layer}$ in powers of $\eps$. We deduce in Section~\ref{A2} the equations satisfied by the asymptotics $\WW_{n}$ \Bk and the 
terms $ \EE_{n}^\pm$. Then we derive explicitly the first terms in Section~\ref{A3}. %

\Bk

\subsection{Expansion of differential operators inside the thin layer}
\label{AH1}

Due to the small thickness of the conductor the derivatives in normal and the tangential directions scale differently in $\eps$.
Hence, it is convenient to use the {\em normal coordinate system}  $\left(y_{\alpha},h\right)$  in $\Otw$. %
For this coordinate system we call 
$D_{\alpha}$ the covariant derivative on the mean surface $\Gamma$ and $\partial_{3}^h$ is  the partial derivative with respect to the normal coordinate $y_{3}=h$. 
Let furthermore $a_{\alpha\beta}(h)$ be the metric tensor of the manifold $\Gamma_{h}$, which is the surface contained in $\Omega^\eps_\layer$ at a distance $h$ of $\Gamma$. 
The metric tensor in such a coordinate system writes~\cite[App. A.1, Eq. (A.7)]{CDFP11}
\begin{equation}
\label{Emett}
a_{\alpha\beta}(h)=a_{\alpha\beta}-2 b_{\alpha\beta} h + b_{\alpha}^{\gamma}b_{\gamma\beta} h^2 \, ,
\end{equation}
 and its inverse expands in power series of $h$
\begin{equation*}
a^{\alpha\beta}(h)=a^{\alpha\beta}+2 b^{\alpha\beta} h + \mathcal{O}( h^2) \, .
\end{equation*}
Here  $b_{\alpha\beta}$ is the curvature tensor on $\Gamma$ and $b_{\alpha}^{\gamma}= a^{\gamma \beta}b_{\beta\alpha}$.
Subsequently, we use a property of the covariant derivative, that it acts on scalar functions $\fke$ like the partial derivative:  $D_{\alpha} \fke=\partial_{\alpha} \fke$.

We denote by $\LL(y_{\alpha}, h;D_{\alpha}, \partial_{3}^h)$ the second order Maxwell operator 
$$
\rot\rot- \kappa^{2}_{\layer} \,\Id
$$ 
in $\Otw$ in the normal coordinate system. We denote by $\D(y_{\alpha}, h;D_{\alpha}, \partial_{3}^h)$ the divergence operator in $\Otw$ \Bk and by $\BB(y_{\alpha} ,h; D_{\alpha}, \partial_{3}^h) = (\BB_\alpha(y_{\alpha} ,h; D_{\alpha}, \partial_{3}^h), 0)$ the tangent trace operator $\rot\cdot\times\nn$ on $\Gamma\cup\Gamma_{\eps}$,
with
\begin{equation*}
\BB_\alpha (y_{\alpha} ,h; D_{\alpha}, \partial_{3}^h) \WW
=\partial_{3}^h \cW_{\alpha} - D_{\alpha} \fke \ ,
\end{equation*}
for $\WW=(\cW_{\alpha},\fke)$, see~\cite[App. A, \S A.4]{CDFP11}. %
The operators $\LL$ and $\BB$ expand in power series of $h$ with intrinsic coefficients with respect to $\Gamma$, see \cite{CDFP11}. 

Now, we scale the normal coordinate $Y_{3}=\eps^{-1} h$ to obtain a coordinate, this is $Y_3$, which does not change with $\eps$.
We use from now on the same {symbol} $\WW$ for three-dimensional one-form field in these scaled coordinates
and we call $\LL[\eps]$, $\BB[\eps]$ and $\D[\eps]$  the respective three-dimensional harmonic Maxwell and divergence operators in $\Otw$. 
These operators expand in powers of $\eps$ 
\begin{equation*}
 \LL[\eps]={\eps}^{-2}\di\sum_{n=0}^{\infty} \eps^n\LL^{n}\ , \quad  \BB[\eps]= {\eps}^{-1}\BB^{0}+\BB^{1} 
 \quad \mbox{and} \quad   \D[\eps]={\eps}^{-1}\di\sum_{n=0}^{\infty} \eps^n\D^{n}  
  \ ,
\end{equation*}
whose coefficients are intrinsic operators on $\Gamma$, which are completely determined by the shape of $\Gamma$ and 
the material parameters of the thin layer. \Bk 
We denote by $L_{\alpha}^n$ and $B_{\alpha}^n$ the surface components of $\LL^{n}$ and $\BB^{n}$.
With the summation convention of repeated two dimensional indices (represented by greek letters), we have (compare with Eq. (A.10) in~\cite[App. A.1]{CDFP11})
\begin{equation}
L_{\alpha}^0(\WW)=-\partial_{3}^2 \cW_{\alpha}
  \ \  \mbox{and}  \ \ 
L_{\alpha}^1(\WW)=-2b_{\alpha}^\beta\partial_{3} \cW_{\beta} 
+ \partial_{3} D_{\alpha} \fke 
+  b_{\beta}^{\beta}\partial_{3}\cW_{\alpha}\ ,
\label{EHsL0L1}
\end{equation}
where  $\partial_{3}$ is the partial derivative with respect to $Y_{3}$, and using~\cite[Eq. (5.39), chap 5]{Pe09} we infer
\begin{equation}
\label{EHsL2}
L^2_\alpha(\WW) = -2 \partial_{3} \left( Y_{3} (b^{2})_{\alpha}^{\beta} \cW_{\beta}\right) + 2\left(D_{\alpha} \gamma_{\beta}^{\beta}(\WW) -D_{\beta} \gamma_{\alpha}^{\beta}(\WW)   \right) -  \kappa^{2}_{\layer} \cW_{\alpha} \ .
\end{equation}
Here $(b^{2})_{\alpha}^{\beta}=b_{\alpha}^{\gamma}b_{\gamma}^{\beta}$, $\gamma_{\alpha\beta}(\WW)=\frac12 (D_{\alpha}\cW_\beta + D_{\beta}\cW_\alpha)-b_{\alpha\beta}\fke$ is 
the change of metric tensor and $\gamma_{\alpha}^{\alpha}= a^{\alpha \beta} \gamma_{\alpha \beta}$.  Furthermore we have~\cite[App. A.1, Eq. (A.28)]{CDFP11} 
\begin{equation}
\label{EB0B1}
B_{\alpha}^0(\WW)= \partial_{3}\cW_{\alpha}   \quad \mbox{and}  \quad B_{\alpha}^1(\WW)= -D_{\alpha} \fke \ ,
\end{equation}
and we have 
\begin{equation}
\label{ED0D1}
\D^0(\WW)= \partial_{3}\fke 
  \quad \mbox{and}  \quad
\D^1(\WW)= \gamma_{\alpha}^{\alpha}(\WW) \ .
\end{equation} 
\Bk

\subsection{Equations for the coefficients of the electric field}
\label{A2}

Writing the partial differential equations~\eqref{eq:EME:2}-\eqref{eq:EME:div} in the thin layer $\Otw$, and the transmission conditions~\eqref{eq:EME:3}  on $\Gamma^{\eps}$ and $\Gamma$
in the semi-scaled local coordinate system $(y_{\alpha},Y_{3})$ we find that 
the asymptotics $\WW_{j}$ \Bk and the terms $ \EE_{j}^\pm$ 
of the electric field 
satisfy the following system \Bk (with $I=(0,1)$)
\begin{subequations}
\label{eq:ELW,EBW}
\begin{align}
\begin{split}
\label{ELW}
&\eps^{-2} \LL[\eps] \sum_{j\geqslant0} \eps^j\WW_j(y_{\alpha},Y_{3}) +(\rot \rot  -\kappa^{2}_{\layer} \Id )\sum_{j\geqslant0} \eps^j \EE^{+}_{j} (\xx)  = 0  \ ,
\\
& \quad
    \mbox{for a.e. } (y_{\alpha},Y_{3})\in \Gamma \times I\ ,  \xx\in\Otw\ ,
    \end{split}
  \\[1.0ex]
   \begin{split}
   \label{EDW}
&\eps^{-1} \D[\eps] \sum_{j\geqslant0} \eps^j\WW_j(y_{\alpha},Y_{3}) +\Div\sum_{j\geqslant0} \eps^j \EE^{+}_{j} (\xx) = 0  \ ,
\\
& \quad
    \mbox{for a.e. } (y_{\alpha},Y_{3})\in \Gamma \times I\ ,  \xx\in\Otw\ ,
        \end{split}
   \\[1.0ex]
  \label{EBW}
& \BB[\eps] \sum_{j\geqslant0} \eps^j\WW_j(y_{\alpha},1) + \sum_{j\geqslant0} \eps^j \rot \EE^{+}_{j} \times\nn =  \frac{\mu_{\layer}}{\mu_{+}}\sum_{j\geqslant0} \eps^j \rot \EE^{+}_{j} \times\nn  \ ,\quad
  \mbox{on} \quad \Gamma^{\eps}\ .
  \\[1.0ex]
   \label{EBW-}
      \begin{split}
&  \frac1{\mu_{\layer}}\BB[\eps] \sum_{j\geqslant0} \eps^j\WW_j(y_{\alpha},0)+   \frac1{\mu_{\layer}} \sum_{j\geqslant0} \eps^j \rot \EE^{+}_{j} \times\nn
\\
& 
=  \frac{1}{\mu_{-}}\sum_{j\geqslant0} \eps^j \rot \EE_{j}^- \times\nn  \ ,
   \quad \mbox{on} \quad \Gamma
  \ .        \end{split}
\end{align}
\end{subequations}
\Bk


In what follows, it is convenient to define $\EE_n$ for $n\in\mathbb{N}$ by $\EE_n=\EE_{n}^+$ in $\Omega_{+}$, and $\EE_n=\EE_{n}^-$ in $\Omega_{-}$. We assume that the terms $\EE_{n}$  are regular in $\Otw$. Then we can use the Taylor expansions of the fields $\EE_{n} : \Otw \longrightarrow  \C^{3}$ and we obtain, for a.e. $\xx\in\Otw$,
\begin{subequations}
\label{eq:Taylor0}
\begin{align} 
\EE_{n} (\xx) &= 
\EE_{n} (\xt) +\eps Y_{3} \partial_h\EE_{n} (\xt) + \cdots \ ,
\\
\rot\rot \EE_{n} (\xx)& = 
\rot\rot\EE_{n} (\xt) +\eps Y_{3} \partial_h\rot\rot\EE_{n} (\xt) + \cdots \ ,
\\
\Div\EE_{n} (\xx)& = 
\Div\EE_{n} (\xt) +\eps Y_{3} \partial_h\Div\EE_{n} (\xt) + \cdots \ ,
\end{align}
\end{subequations}
where $\xt=P(\xx)$ is the orthogonal projection of $\xx$ on $\Gamma$, $Y_{3}=\dfrac{h}{\eps}$ where $h=\mathrm{dist}(\xx,\Gamma)$
.   We also have the following  Taylor expansions
\begin{subequations}
\label{eq:Taylor}
\begin{align} 
\EE_{n} \times\nn |_{h=\eps} &= \hspace{1.7em}
\EE_{n}  \times\nn |_{0^+}+\eps\partial_{h}\EE_{n} \times\nn |_{0^+}+\cdots\ ,
\\
\EE_{n} \cdot\nn |_{h=\eps} &= \hspace{1.7em}
\EE_{n}  \cdot\nn |_{0^+}+\eps\partial_{h}\EE_{n} \cdot\nn |_{0^+}+\cdots\ ,
\\
\rot \EE_{n} \times\nn |_{h=\eps} &= \rot 
\EE_{n}\times\nn |_{0^+}+\eps \partial_{h}  \rot \EE_{n} \times\nn |_{0^+} +\cdots\ ,
\end{align}
\end{subequations}
where $\cdot|_{0^+}$ means the limit for positive $h \to 0$.

Then according to the systems~\eqref{eq:EME}-\eqref{eq:ELW,EBW} and the transmission conditions~\eqref{eq:EMEbis}  on $\Gamma$ and  $\Gamma^{\eps}$,
using the expressions of the operators $\LL^0$, $\BB^0$, $\BB^1$, $\D^{0}$, and using the Taylor expansions~\eqref{eq:Taylor}-\eqref{eq:Taylor}, we deduce that the  terms $\WW_n=(\cW_n,\fke_n )$ \Bk and $\EE_{n}$ have to satisfy, for all  $n \geq 0$ 
\Bk

\begin{subequations}\label{EE:n}
\begin{align*}
  %
& -\partial_{3}^2 \cW_{n,\alpha}   =  -\sum_{j=1}^n L^j_\alpha(\WW_{n-j}) 
 - \di\sum_{j=2}^{n} \frac1{\fact{(j-2)}}\, Y^{j-2}_{3}\partial_{h}^{j-2}\left( (\rot\rot  - \kappa^{2}_{\layer}\Id ) \EE^+_{n-j}\right)_{\alpha}(\xt)
 	\quad \mbox{in}\quad \Gamma \times I,
	\\[0.8ex] 
&	 \partial_{3}\fke_n=  -\sum_{j=1}^n \D^j(\WW_{n-j})  - \di\sum_{j=1}^{n} \frac1{\fact{(j-1)}}\, Y^{j-1}_{3}\partial_{h}^{j-1} \left(\Div\EE^+_{n-j}\right)(\xt)\quad \mbox{in}\quad \Gamma \times I,
  \\[0.8ex]  %
& \partial_{3}\cW_{n,\alpha}|_{1  } = D_\alpha\fke_{n-1}|_{1  } 
 + (\frac{\mu_{\layer}}{\mu_+}-1)\sum_{j=1}^{n}\frac{1}{\fact{(j-1)}}\partial_{h}^{j-1}\hspace{-0.25em}\left(\curl \EE^+_{n-j}\times\nn\right)_\alpha  \mbox{on} \quad \Gamma,
  \\[0.8ex] 
& \partial_{3}\cW_{n,\alpha}|_{0} = D_\alpha\fke_{n-1}|_{0} 
 + \frac{\mu_{\layer}}{\mu_-} 
 \left(\curl \EE^-_{n-1}\times\nn\right)_\alpha-  \left(\curl \EE^+_{n-1}\times\nn\right)_\alpha  \mbox{on} \quad \Gamma,
  \\[0.8ex] 
  %
 &  \rot\rot  \EE^{-}_{n} - \kappa_{-}^2 \EE^{-}_{n} = 0
\quad    \mbox{in}\quad \Omega_{-},
\\[0.8ex]
 &  \rot\rot  \EE^{+}_{n} - \kappa_{+}^2 \EE^{+}_{n} = \delta^0_n i\omega\mu_+ \JJ
\quad    \mbox{in}\quad \Omega_{+},
\\[0.8ex]
 %
&\EE_{n}^{-} \times\nn   = \EE_{n}^{+} \times\nn  +  \WW_{n}\times\nn|_{0}  \quad\mbox{on}\quad \Gamma,
\\[0.8ex]
 %
& \WW_{n}\times\nn|_{1}=0  \quad\mbox{on}\quad \Gamma,
 \\[0.8ex] 
& \left(i \omega \epsilon_{-}-\sigma_{-} \right) \EE_{n}^{-}\cdot \nn   
=  \left(i \omega \epsilon_{\layer}-\sigma_{\layer} \right) \left(\EE_{n}^{+} \cdot\nn  +  \fke_{n}|_{0}  \right)\quad\mbox{on}\quad \Gamma,
   \\[0.8ex] 
    \begin{split}
& \left(i \omega \epsilon_{\layer}-\sigma_{\layer} \right) \left(\EE_{n}^{+}\cdot \nn   + \di\sum_{j=1}^{n} \frac1{\fact{j}} \partial_{h}^{j}\EE_{n-j}^{+}\cdot\nn+  \fke_{n}|_{1}   \right)
\\
&=  \left(i \omega \epsilon_{+}-\sigma_{+} \right) \left(\EE_{n}^{+} \cdot\nn+\di\sum_{j=1}^{n} \frac1{\fact{j}} \partial_{h}^{j}\EE_{n-j}^{+}\cdot\nn \right) \quad\mbox{on}\quad \Gamma,
 \end{split}
\\[0.8ex]
  %
&\curl{\EE_{n}^+} \times \nn - i \kappa_{+}\, \nn \times \EE_{n}^+ \times \nn = 0 \quad \mbox{on} \quad \partial\Om \ , 
\end{align*}
\end{subequations}
where $\cdot|_{1}$ and  $\cdot|_{0}$  abbreviate  the traces on $Y_3  =1$ and $Y_{3}=0$,  respectively. 
In the previous system we use the convention that the sums are $0$ when $n=0$. 
\begin{rem}
\label{R1}
Note that the condition $\WW_{n}\times\nn|_{1}=0$  on $\Gamma$, comes from the transmission condition \eqref{eq:EME:4} set on $\Gamma^{\eps}$ since the tangential traces $\EE_{n}^{+} \times\nn$ are continuous across $\Gamma^{\eps}$.
\end{rem}

\subsection{First terms of the asymptotics}
\label{A3}
In the previous section we have derived the coupled systems for the terms of the WKB expansion to any order~n. Hence we can determine now the first terms $\WW_n=(\cW_n,\fke_n )$ and $\EE_{n}$ by induction.

\subsubsection*{The coupled system of order 0} %
For $n=0$ in the previous system, we find that  $\WW_0=(\cW_0,\fke_0 )$ and $\EE_0$ satisfy
\begin{subequations}\label{EE:0}
\begin{align}
\label{EE:0:1}
  -\partial_{3}^2 \cW_{0,\alpha}
   &=  0  && \mbox{in}\quad \Gamma \times I\ ,
   \\[0.8ex]
  \label{EE:0:2}
   \partial_{3}\fke_0 &=  0
 &&\mbox{in}\quad \Gamma \times I\ ,
 \\[0.8ex]
  \label{EE:0:3}
  \partial_{3}\cW_{0,\alpha}|_{1}=  \partial_{3}\cW_{0,\alpha}|_{0} &= 0  &&\mbox{on} \quad \Gamma\ ,
  \\[0.8ex]
  \label{EE:0:4}
   \rot\rot  \EE^{-}_{0} - \kappa_{-}^2 \EE^{-}_{0} &= 0   &&\mbox{in}\quad \Omega_{-}\ ,
    \\[0.8ex]
  \label{EE:0:4bis}
   \rot\rot  \EE^{+}_{0} - \kappa_{+}^2 \EE^{+}_{0} &= i \omega   \mu_{+}\JJ  &&\mbox{in}\quad \Omega_{+}\ ,
\\[0.8ex]
 \label{EE:0:5}
\EE_{0}^{-} \times\nn& =\EE_{0}^{+} \times\nn  + \WW_{0}\times\nn|_{0} && \mbox{on}\quad \Gamma\ ,
\\[0.8ex]
 \label{EE:0:7}
\WW_{0}\times\nn|_{1} &=0 && \mbox{on}\quad \Gamma\ ,
\\[0.8ex]
   \label{EE:0:8}
 \left(i \omega \epsilon_{-}-\sigma_{-} \right) \EE_{0}^{-}\cdot \nn &  =  \left(i \omega \epsilon_{\layer}-\sigma_{\layer} \right)\left( \EE_{0}^{+} \cdot\nn  +  \fke_{0}|_{0}  \right)&&  \mbox{on}\quad \Gamma \ ,
  \\[0.8ex]
   \label{EE:0:9}
 \left(i \omega \epsilon_{\layer}-\sigma_{\layer} \right) \left(\EE_{0}^{+}\cdot \nn +  \fke_{0}|_{1}   \right)&=  \left(i \omega \epsilon_{+}-\sigma_{+} \right) \EE_{0}^{+} \cdot\nn &&  \mbox{on}\quad \Gamma \ ,
\\[0.8ex]
  \label{EE:0:6}
\curl{\EE_{0}^+} \times \nn - i \kappa_{+}\, \nn \times \EE_{0}^+ \times \nn &= 0 && \mbox{on} \quad \partial\Om\ .
\end{align}
\end{subequations}
Obviously,~\eqref{EE:0:1},~\eqref{EE:0:7} and~\eqref{EE:0:3} imply that   $\cW_{0,\alpha} = 0$ in $\Gamma \times I$. We deduce from~\eqref{EE:0:5} that  
\begin{equation}
\label{EE05}
\EE_{0}^{-} \times\nn =\EE_{0}^{+} \times\nn  \quad  \mbox{on}\quad \Gamma \ .
\end{equation}
In view of~\eqref{EE:0:2} we can assert that 
$\fke_0$ does not depend on the variable $Y_{3}$ in $\Gamma \times I$. We infer from~\eqref{EE:0:8} and~\eqref{EE:0:9} that 
\begin{equation}
\label{EE07}
\left(i \omega \epsilon_{-}-\sigma_{-} \right) \EE_{0}^{-}\cdot \nn =  \left(i \omega \epsilon_{+}-\sigma_{+} \right) \EE_{0}^{+} \cdot\nn  \quad \mbox{on}\quad \Gamma \ ,
\end{equation}
 and we deduce that 
\begin{equation}
\label{EE08}
 \fke_{0}(y_{\alpha})=\left(\dfrac{i \omega\epsilon_{+}-\sigma_{+}}{i \omega  \epsilon_{\layer}-\sigma_{\layer}}-1\right)\EE_{0}^{+}\cdot \nn |_{\Gamma^{+}}
 =\left(\dfrac{ i \omega  \epsilon_{-}-\sigma_{-}}{i \omega \epsilon_{\layer}-\sigma_{\layer}}-\dfrac{ i \omega  \epsilon_{-}-\sigma_{-}}{i \omega \epsilon_{+}-\sigma_{+}}\right)\EE_{0}^{-}\cdot \nn |_{\Gamma^{-}}\ , 
\end{equation}
where  $\EE_{0}$ has to be determined.  
  
  %

\subsubsection*{The coupled system of order 1} 
Then   in the same way as above we find that $\WW_1=(\cW_1,\fke_1 )$ and  $\EE_1$ satisfy \

\begin{subequations}
\label{EE:1}
\begin{align}
  %
\label{EE:1:1}
 -\partial_{3}^2 \cW_{1,\alpha}   =  -L^1_\alpha(\WW_{0})  	\quad \mbox{in}\quad \Gamma \times I,
	\\
	\label{EE:1:0}
	 \partial_{3}\fke_1=  - \D^1(\WW_{0})  - \Div\EE^+_{0}(\xt)\quad \mbox{in}\quad \Gamma \times I,
  \\[0.8ex]  %
  \label{EE:1:2}
 \partial_{3}\cW_{1,\alpha}|_{1  } = D_\alpha\fke_{0}|_{1  } 
 + \left(\frac{\mu_{\layer}}{\mu_+}-1\right)\left(\curl \EE^+_{0}\times\nn\right)_\alpha  \mbox{on} \quad \Gamma,
  \\[0.8ex]  %
  \label{EE:1:3}
 \partial_{3}\cW_{1,\alpha}|_{0} = D_\alpha\fke_{0}|_{0} 
 + \frac{\mu_{\layer}}{\mu_-} 
 \left(\curl \EE^-_{0}\times\nn\right)_\alpha-  \left(\curl \EE^+_{0}\times\nn\right)_\alpha  \mbox{on} \quad \Gamma,
  \\[0.8ex]
  %
  \label{EE:1:6}
   \rot\rot  \EE^{\pm}_{1} - \kappa_{\pm}^2 \EE^{\pm}_{1} = 0
\quad    \mbox{in}\quad \Omega_{\pm},
\\[0.8ex]
 \label{EE:1:4}
\EE_{1}^{-} \times\nn   = \EE_{1}^{+} \times\nn  +  \WW_{1}\times\nn|_{0}  \quad\mbox{on}\quad \Gamma,
\\[0.8ex]
 \label{EE:1:5}
 \WW_{1}\times\nn|_{1}=0  \quad\mbox{on}\quad \Gamma,
 \\[0.8ex]
  \label{EE:1:7}
 \left(i \omega \epsilon_{-}-\sigma_{-} \right) \EE_{1}^{-}\cdot \nn   =  \left(i \omega \epsilon_{\layer}-\sigma_{\layer} \right) \left(\EE_{1}^{+} \cdot\nn  +  \fke_{1}|_{0}  \right)\quad\mbox{on}\quad \Gamma,
 \\[0.8ex]
   \label{EE:1:8}
 \left(i \omega \epsilon_{\layer}-\sigma_{\layer} \right) \left(\EE_{1}^{+}\cdot \nn   + \partial_{h}\EE_{0}^{+}\cdot\nn+  \fke_{1}|_{1}   \right)=  \left(i \omega \epsilon_{+}-\sigma_{+} \right) \left(\EE_{1}^{+} \cdot\nn+\partial_{h}\EE_{0}^{+}\cdot\nn \right) \quad\mbox{on}\quad \Gamma,
\\[0.8ex]
  \label{EE:1:9}
\curl{\EE_{1}^+} \times \nn - i \kappa_{+}\, \nn \times \EE_{1}^+ \times \nn = 0 \quad \mbox{on} \quad \partial\Om.
\end{align}
\end{subequations}

The right hand side of \eqref{EE:1:1} is zero since  $\cW_{0,\alpha} = 0$ in $\Gamma \times I$ and $\fke_0$ depends only  on the variables $y_{\alpha}$. Hence $\partial_{3}\cW_{1,\alpha}$ depends only  on the variables $y_{\alpha}$, and we infer that the tangential components  of $\cW_{1,\alpha}$ are given for $Y_{3}\in I $ by 
\begin{equation}
\label{EcW1}
\cW_{1,\alpha}(\cdot, Y_{3})=a_{\alpha} +  b_{\alpha}Y_{3}\ , 
\end{equation}
where $a_{\alpha}$ and $b_{\alpha}$ depend only on the variables $y_{\beta}$, and have to be determined. 
Furthermore since $D_\alpha\fke_{0}|_{0}= D_\alpha\fke_{0}|_{1}$ we deduce from \eqref{EE:1:2}-\eqref{EE:1:3}  that 
\begin{equation*}
 \left(\frac{\mu_{\layer}}{\mu_+}-1\right)\left(\curl \EE^+_{0}\times\nn\right)_\alpha  = \frac{\mu_{\layer}}{\mu_-} 
 \left(\curl \EE^-_{0}\times\nn\right)_\alpha-  \left(\curl \EE^+_{0}\times\nn\right)_\alpha  \mbox{on} \quad \Gamma .
\end{equation*}
We infer the following transmission conditions of $\EE_{0}$ across  $\Gamma$ 
 \begin{equation}
 \label{EE09}
\frac1{\mu_+}\left(\curl \EE^+_{0}\times\nn\right)  = \frac{1}{\mu_-} 
 \left(\curl \EE^-_{0}\times\nn\right)  \mbox{on} \quad \Gamma \ ,
\end{equation}
and using \eqref{EE:0:4}-\eqref{EE:0:4bis}-\eqref{EE05}-\eqref{EE07}-\eqref{EE:0:6} we conclude that the term $\EE_0$ solves Problem \eqref{0E0}. Therefore the terms of order $0$ are entirely determined.

\vspace{2mm}
Then using \eqref{EE:1:5}, we deduce successively that we have  $a_{\alpha}+b_{\alpha}=0$ and 
\begin{equation}
\label{EEcW1}
\cW_{1,\alpha}(\cdot, Y_{3})=a_{\alpha} (1- Y_{3}) \ ,
\end{equation}
 where $a_{\alpha}= -\partial_{3}\cW_{1,\alpha}|_{1} =- \partial_{3}\cW_{1,\alpha}|_{0}$ can be determined by using  \eqref{EE:1:2} or \eqref{EE:1:3}. Using   \eqref{EE:1:2}, we infer
\begin{equation}
\label{aalpha}
a_{\alpha}= - D_\alpha\fke_{0}|_{1}  - \left(\frac{\mu_{\layer}}{\mu_+}-1\right)\left(\curl \EE^+_{0}\times\nn\right)_\alpha \quad\mbox{on}\quad \Gamma \ . 
\end{equation}
Then, using  \eqref{EE:1:4}, we have 
\begin{equation*}
 \EE_{1}^{+} \times\nn  = \EE_{1}^{-} \times\nn  -  \WW_{1}\times\nn|_{0}  \quad\mbox{on}\quad \Gamma \ .
\end{equation*}
We infer successively  (using \eqref{EE09}) 
\begin{equation*}
 \EE_{1}^{+} \times\nn  = \EE_{1}^{-} \times\nn + \nabla_{\Gamma} \fke_{0}|_{1} \times\nn
 - (\mu_{\layer}-\mu_{+}) \left(\frac1{\mu_{-}}\curl \EE^{-}_{0}\right) _{\T} 
   \quad\mbox{on}\quad \Gamma \ ,
\end{equation*}
 and 
\begin{multline}
\label{**}
 \EE_{1}^{+} \times\nn  = \EE_{1}^{-} \times\nn
 \\
  + \left(\frac{\mu_{\layer}}{\kappa_{\layer}^2}-\frac{\mu_{+}}{\kappa_{+}^2}\right)\veccurl_\Gamma \curl_{\Gamma} (\frac1{\mu_{-}}\curl \EE^{-}_{0})_{\T}
- (\mu_{\layer}-\mu_{+}) \left(\frac1{\mu_{-}}\curl \EE^{-}_{0}\right) _{\T} 
   \quad\mbox{on}\quad \Gamma \ ,
\end{multline}
by using \eqref{EE08} together with the following identity 
\begin{equation*}
\frac{\kappa_{-}^2}{\mu_{-}}\nabla_\Gamma\left(\EE_0^- \cdot \nn\right)=
\nn\times \veccurl_\Gamma\curl_\Gamma \left(\frac{1}{\mu_{-}}\curl\EE_0^-\right)_\T  \quad\mbox{on}\quad \Gamma \ ,
\end{equation*}
which comes from \eqref{EE:0:4} (see, {\it e.g.},~\cite[Section 3.2]{Durufle2014}). 

According to \eqref{EE:1:0},  $\partial_{3}\fke_1$ depends only on the variables $y_{\alpha}$ 
since $\WW_0$ depends only on the variables $y_{\alpha}$.  We deduce that $\fke_1$ is given for $Y_{3}\in I $ by 
\begin{equation}
\label{Efke1}
\fke_1(\cdot, Y_{3})=a +  b Y_{3}\ , 
\end{equation}
where functions $a$ and $b$  depend only on the variables $y_{\beta}$, and have to be determined.  Then using  \eqref{EE:1:7}-\eqref{EE:1:8}, since $\fke_1|_{0} =a$ and $\fke_1|_{1} =a +  b$, we infer successively 
\begin{equation}
\label{Efke1a}
a=\dfrac{i \omega\epsilon_{-}-\sigma_{-}}{i \omega  \epsilon_{\layer}-\sigma_{\layer}} \, \EE_{1}^{-}\cdot \nn |_{\Gamma} - \EE_{1}^{+}\cdot \nn |_{\Gamma} 
\ ,
\end{equation}
\begin{equation*}
a+b= \left(\dfrac{i \omega\epsilon_{+}-\sigma_{+}}{i \omega  \epsilon_{\layer}-\sigma_{\layer}}-1\right)\left(\EE_{1}^{+}\cdot \nn |_{\Gamma} + \partial_{h}\EE_{0}^{+}\cdot \nn |_{\Gamma}\right)\ ,
\end{equation*}
and
\begin{equation}
\label{Efke1b}
b=  \left(\dfrac{i \omega\epsilon_{+}-\sigma_{+}}{i \omega  \epsilon_{\layer}-\sigma_{\layer}}\right)\EE_{1}^{+}\cdot \nn |_{\Gamma} - \left(\dfrac{i \omega\epsilon_{-}-\sigma_{-}}{i \omega  \epsilon_{\layer}-\sigma_{\layer}} -1 \right)\EE_{1}^{-}\cdot \nn |_{\Gamma}
+ \left(\dfrac{i \omega\epsilon_{+}-\sigma_{+}}{i \omega  \epsilon_{\layer}-\sigma_{\layer}}-1\right) \partial_{h}\EE_{0}^{+}\cdot \nn |_{\Gamma}\ ,
\end{equation}
where $\EE_{1}$ has to be determined. 
\Bk

\subsubsection*{The coupled system of order 2} 
In the same way as above we find that the asymptotics $\WW_2=(\cW_2,\fke_2)$ and  $\EE_2$ satisfy 

\begin{subequations}\label{EE:2}
\begin{align}
\label{EE:2:1}
 &-\partial_{3}^2 \cW_{2,\alpha}   =  - L^1_\alpha(\WW_{1}) -L^2_\alpha(\WW_{0}) 
 - \left( (\rot\rot  -\kappa^{2}_{\layer} \Id ) \EE^+_{0}\right)_{\alpha}(\xt)
 	\quad \mbox{in}\quad \Gamma \times I,
	 \\[0.8ex]
&	 \partial_{3}\fke_2=  -\sum_{j=1}^2 \D^j(\WW_{2-j})  - \di\sum_{j=1}^{2}  Y^{j-1}_{3}\partial_{h}^{j-1} \left(\Div\EE^+_{2-j}\right)(\xt)\quad \mbox{in}\quad \Gamma \times I,
  \\[0.8ex]  %
  \label{EE:2:2}
& \partial_{3}\cW_{2,\alpha}|_{1  } = D_\alpha\fke_{1}|_{1  } 
 + (\frac{\mu_{\layer}}{\mu_+}-1)\sum_{j=1}^{2}\partial_{h}^{j-1}\hspace{-0.25em}\left(\curl \EE^+_{2-j}\times\nn\right)_\alpha  \mbox{on} \quad \Gamma,
  \\[0.8ex]  %
  \label{EE:2:3}
& \partial_{3}\cW_{2,\alpha}|_{0} = D_\alpha\fke_{1}|_{0} 
 + \frac{\mu_{\layer}}{\mu_-} 
 \left(\curl \EE^-_{1}\times\nn\right)_\alpha-  \left(\curl \EE^+_{1}\times\nn\right)_\alpha  \mbox{on} \quad \Gamma,
  \\[0.8ex]
  %
 &  \rot\rot  \EE^{\pm}_{2} - \kappa_{\pm}^2 \EE^{\pm}_{2} =0
\quad    \mbox{in}\quad \Omega_{\pm},
\\[0.8ex]
 %
&\EE_{2}^{-} \times\nn   = \EE_{2}^{+} \times\nn  +  \WW_{2}\times\nn|_{0}  \quad\mbox{on}\quad \Gamma,
\\[0.8ex]
 %
& \WW_{2}\times\nn|_{1}=0  \quad\mbox{on}\quad \Gamma,
 \\[0.8ex]
& \left(i \omega \epsilon_{-}-\sigma_{-} \right) \EE_{2}^{-}\cdot \nn   =  \left(i \omega \epsilon_{\layer}-\sigma_{\layer} \right) \left(\EE_{2}^{+} \cdot\nn  +  \fke_{2}|_{0}  \right)\quad\mbox{on}\quad \Gamma,
   \\[0.8ex]
 \begin{split}
 & \left(i \omega \epsilon_{\layer}-\sigma_{\layer} \right) \left(\EE_{2}^{+}\cdot \nn   + \di\sum_{j=1}^{2} \frac1{\fact{j}} \partial_{h}^{j}\EE_{2-j}^{+}\cdot\nn+  \fke_{2}|_{1}   \right)
  \\[0.8ex]
& =  \left(i \omega \epsilon_{+}-\sigma_{+} \right) \left(\EE_{2}^{+} \cdot\nn+\di\sum_{j=1}^{2} \frac1{\fact{j}} \partial_{h}^{j}\EE_{2-j}^{+}\cdot\nn \right) \quad\mbox{on}\quad \Gamma,
 \end{split}
\\
  %
& \curl{\EE_{2}^+} \times \nn - i \kappa_{+}\, \nn \times \EE_{2}^+ \times \nn = 0 \quad \mbox{on} \quad \partial\Om.
\end{align}
\end{subequations}
First, by definition \eqref{EHsL2}, the term $L^2_\alpha(\WW_{0})$ which appears in the right-hand side of identity \eqref{EE:2:1}
 is independent of the variable $Y_{3}$  since  $\cW_0=0$ in $ \Gamma \times I$ and $\fke_{0}$  is independent of the variable $Y_{3}$. Second, according to \eqref{EHsL0L1}-\eqref{EcW1}-\eqref{Efke1}, the term $L^1_\alpha(\WW_{1})$ is also independent of the variable $Y_{3}$  in $ \Gamma \times I$. Then in view of \eqref{EE:2:1} we deduce that the term $\partial_{3}^2 \cW_{2,\alpha}$ is independent of the variable $Y_{3}$ in $ \Gamma \times I$. Hence we have $ \partial_{3}\cW_{2,\alpha}|_1 = \partial_{3}\cW_{2,\alpha}|_0$ and  \eqref{EE:2:2}-\eqref{EE:2:3} yield to the following transmission conditions 
\begin{multline*}
\frac{1}{\mu_+}\left(\curl \EE^+_{1}\times\nn\right)-\frac{1}{\mu_-}\left(\curl \EE^-_{1}\times\nn\right)=\frac1{\mu_{\layer}}\left(\nabla_{\Gamma}\fke_{1}|_{0}- \nabla_{\Gamma}\fke_{1}|_{1}\right)
\\
+ (\frac1{\mu_{\layer}}- \frac1{\mu_+})\partial_{h}\hspace{-0.25em}\left(\curl \EE^+_{0}\times\nn\right) \quad \mbox{on} \quad  \Gamma \ .
\end{multline*}
We infer the transmission conditions \eqref{TC2} by proving hereafter the following two identities 
\begin{equation}
\label{Eputain}
 \frac{1}{\mu_{\layer}} \left( \nabla_{\Gamma}\fke_{1}|_{0}- \nabla_{\Gamma}\fke_{1}|_{1}\right)  = 
 \frac{2}{\mu_{\layer}} \frac{\mu_{+}}{\kappa_{+}^2} \left(1-\dfrac{i \omega\epsilon_{+}-\sigma_{+}}{i \omega  \epsilon_{\layer}-\sigma_{\layer}}\right) 
\nn\times\veccurl_{\Gamma} \left( \mathcal{H}\curl_{\Gamma} (\frac1{\mu_{+}}\curl \EE_{0}^{+})_{\T}\right) \ ,
\end{equation}
and
\begin{equation}
\label{Putain}
\partial_{h}\hspace{-0.25em}\left(\curl \EE^+_{0}\times\nn\right) = 
\kappa_{+}^2  \   (\EE^+_{0})_{\T}
+
\veccurl_\Gamma\rot_\Gamma  (\EE^+_{0})_{\T} 
  \quad \mbox{on} \quad \Gamma \ .
\end{equation}
 Finally, in view of \eqref{EE:1:6}-\eqref{**}-\eqref{EE:1:9}, we deduce that $\EE_{1}$ solves Problem \eqref{0E1}.  Furthermore the term $\fke_1$ is determined by \eqref{Efke1}-\eqref{Efke1a}-\eqref{Efke1b}. Therefore the terms of order $1$ are entirely determined.

\Bk 
\begin{proof}[Proof of identities \eqref{Eputain}-\eqref{Putain}]
According to \eqref{EE:1:0}, we have $\partial_{3}\fke_1=-\D^{1}(\WW_{0})$ in $ \Gamma \times I$ since we have $\Div\EE^+_{0}=0$ in a neighborhood of $\Gamma$ by assumption on $\JJ$. Then according to \eqref{ED0D1}, we infer successively 
\begin{equation*}
\partial_{3}\fke_1= -\gamma_{\alpha}^{\alpha}(\WW_{0}) \quad \mbox{and} \quad \partial_{3}\fke_1=b_{\alpha}^{\alpha} \fke_{0}  \quad \mbox{in} \quad    \Gamma \times I
 \end{equation*}
since $\cW_{0}=0$. Using \eqref{EE08} we deduce the following identity 
\begin{equation}
\label{Ee1}
\partial_{3}\fke_1 =  2\mathcal{H} \left(\dfrac{i \omega\epsilon_{+}-\sigma_{+}}{i \omega  \epsilon_{\layer}-\sigma_{\layer}}-1\right)\EE_{0}^{+}\cdot \nn |_{\Gamma} \quad \mbox{in} \quad \Gamma\times I \ 
\end{equation}
since  $b_{\alpha}^{\alpha}  = 2\mathcal{H}$. Then  taking the tangential gradient on $\Gamma$ of each term in \eqref{Ee1} and integrating each term in the variable $Y_{3}$ along $I$, we obtain   
\begin{equation}
\label{putain}
\nabla_{\Gamma}\fke_{1}|_{1} - \nabla_{\Gamma}\fke_{1}|_{0} =  2\left(\dfrac{i \omega\epsilon_{+}-\sigma_{+}}{i \omega  \epsilon_{\layer}-\sigma_{\layer}}-1\right)\nabla_{\Gamma}  \left( \mathcal{H} \EE_{0}^{+}\cdot \nn |_{\Gamma} \right)\quad \mbox{in} \quad \Gamma\times I \ .
\end{equation}
Finally, multiplying each term in \eqref{EE:0:4bis} by $\mathcal{H}$ and then taking successively the normal trace on $\Gamma$ of each term and the tangential gradient on $\Gamma$ of each term, we deduce the following identity 
\begin{equation*}
\kappa_{+}^2\nabla_\Gamma\left( \mathcal{H} \EE_0^+ \cdot \nn\right)=
\nn\times \veccurl_\Gamma\left( \mathcal{H} \curl_\Gamma \left(\curl\EE_0^+\right)_\T \right) \quad\mbox{on}\quad \Gamma \ ,
\end{equation*}
from which we infer identity \eqref{Eputain} by using  \eqref{putain}.

Let us prove now identity \eqref{Putain} 
\begin{equation*}
\partial_{h}\hspace{-0.25em}\left(\curl \EE^+_{0}\times\nn\right) = 
\kappa_{+}^2  \   (\EE^+_{0})_{\T}
+
\veccurl_\Gamma\rot_\Gamma  (\EE^+_{0})_{\T} 
  \quad \mbox{on} \quad \Gamma \ .
\end{equation*}
Following Section~\ref{sec:Asymp}, it is also possible to derive at least formally a WKB expansion for the magnetic field $\HH^\eps$ as follow
\begin{subequations}
\label{eq:MHbis}
\begin{align}
\HH^\eps(\xx)
\approx{\HH}^{0}(\xx)+\eps{\HH}^{1}(\xx)+\eps^{2}{\HH}^{2} (\xx)+\cdots ,  \quad 
\qquad \text{for a.e.}\quad \xx\in \Omega_{-}\cup\Omega^\eps_{+}\ ,
  \\
\HH^\eps(\xx)
\approx {\HH}^{0}(\xx)+\VV_0 \left(y_{\alpha},\frac{h}{\eps}\right) + \eps \left({\HH}^{1}(\xx)+ \VV_1 \left(y_{\alpha},\frac{h}{\eps}\right)  \right) 
+\cdots, 
 \qquad \text{for a.e.}\quad \xx \in \Otw\ .
\end{align}
\end{subequations}
 Then, using Faraday's law, we deduce that  the first terms $\EE^+_{0},\HH^+_{0}$ solve
 \begin{equation*}
\curl \EE^+_{0} =  i \omega \mu_{+}\HH^+_{0} \quad \mbox{in} \quad \Omega_{+} \ .
\end{equation*}
We infer
 \begin{equation}
 \label{e1}
\partial_{h}\hspace{-0.25em}\left(\curl \EE^+_{0}\times\nn\right) =  i \omega \mu_{+}\partial_{h}\hspace{-0.25em}\left(\HH^+_{0}\times\nn\right) 
  \quad \mbox{on} \quad \Gamma \ .
\end{equation}
Furthermore, according to \cite[Appendix A]{Pe16}, we have
\begin{equation}
 \label{e2}
\partial_{h}\hspace{-0.25em}\left(\HH^+_{0}\times\nn\right)  = \left(\curl \HH^+_{0}\right)_{\T}+ \nabla_{\Gamma}  \left(\HH_{0}^{+}\cdot \nn \right)\times\nn  \quad \mbox{on} \quad \Gamma \ .
\end{equation}
Then, using Amp\`ere's law, we check that  the first terms $\EE^+_{0},\HH^+_{0}$ satisfy
\begin{equation}
 \label{e3}
 \left(\curl \HH^+_{0}\right)_{\T} =- \left(i \omega \epsilon_{+}-\sigma_{+} \right) ( \EE^+_{0} )_{\T}   \quad \mbox{on} \quad \Gamma \ .
\end{equation}
Furthermore, using Faraday's and  Amp\`ere's laws, we have
\begin{equation*}
\curl\curl\HH_0^+ -\kappa^2_{+}\HH_0^+ = \curl \JJ  \quad \mbox{in} \quad\Omega_{+} \ , 
\end{equation*}
from which we obtain
 \begin{equation}
  \label{e4}
\kappa_{+}^2\nabla_\Gamma\left(\HH_0^+ \cdot \nn\right)=
\nn\times \veccurl_\Gamma\curl_\Gamma \left(\curl\HH_0^+\right)_\T  \quad\mbox{on}\quad \Gamma \ .
\end{equation}
Finally  we infer  identity \eqref{Putain} by combining the previous identities \eqref{e1}-\eqref{e2}-\eqref{e3}-\eqref{e4}. 
\end{proof}
\begin{rem}
Note also that according to \eqref{Efke1}-\eqref{Efke1a}-\eqref{Ee1} we infer 
\begin{equation*}
\fke_1(\cdot, Y_{3}) =
\dfrac{i \omega\epsilon_{-}-\sigma_{-}}{i \omega  \epsilon_{\layer}-\sigma_{\layer}} \, \EE_{1}^{-}\cdot \nn |_{\Gamma} - \EE_{1}^{+}\cdot \nn |_{\Gamma} +  2\mathcal{H} \left(\dfrac{i \omega\epsilon_{+}-\sigma_{+}}{i \omega  \epsilon_{\layer}-\sigma_{\layer}}-1\right)\EE_{0}^{+}\cdot \nn |_{\Gamma}   Y_{3} \quad \mbox{in}\quad \Gamma \times I \ .
\end{equation*}
We deduce the expression of the term $\Em^1$ in \eqref{EEm1}  by using \eqref{EEcW1}-\eqref{aalpha}.
\end{rem}

 
\subsection{A remark on remainders of the WKB expansion in comparison with the multiscale expansion derived  in~\cite{DPP11}}
\label{S4.4}

The validation of the WKB expansion \eqref{eq:ME} consists in proving estimates for remainders $ \RR^{m;\,\eps}$ defined as 
\begin{subequations}
\label{eq:MER}
\begin{align}
 \RR_{-}^{m;\,\eps} & =
\EE_{-}^\eps- \di\sum_{j =0}^{m} \eps^{j} \EE_{-}^{j} 
,  \quad \text{in}\quad \Omega_{-}  \ ,  \\
\label{eq:ME2R}
 \RR_{\layer}^{m;\,\eps} &
 =\EE_{\layer}^\eps -  \di\sum_{j =0}^{m} \eps^{j} \left(\EE_{+}^{j}+\EE_{\layer}^j\left(\xt,\frac{x_3}{\eps}\right) \right)
, 
\quad \text{in}\quad   \Omega^\eps_{\layer} \ ,  \\ 
 \RR_{+}^{m;\,\eps} &
 =\EE_{+}^\eps - \di\sum_{j =0}^{m} \eps^{j} \EE_{+}^{j} 
,  \quad \text{in}\quad  \Omega_{+}^\eps
 \ .
\end{align}
\end{subequations}
Error estimates have been proved for the multiscale expansion derived in~\cite{DPP11}, see~\cite[Th. 2.9]{DPP11}. By construction of the WKB expansion, the remainder $\RR^{m;\,\eps}$ defines an element of $\bH(\curl,\Omega)$ for smooth data 
since we have 
\begin{equation*}
 \RR_{-}^{m;\,\eps} \times \nn =  \RR_{\layer}^{m;\,\eps} \times \nn \quad \mbox{on} \quad \Gamma \ ,  \qquad \mbox{and} \quad  \RR_{+}^{m;\,\eps} \times \nn= \RR_{\layer}^{m;\,\eps} \times \nn \quad \mbox{on} \quad \Gamma^{\eps} \ ,
\end{equation*}
 see the last system of equations in Section~\ref{A2} and Remark~\ref{R1}. 
This is a main difference with the  multiscale expansion derived in~\cite{DPP11} for which tangential  traces of the remainders  $\tilde\RR^{m;\,\eps}$  are discontinuous across the boundary $\Gamma^{\eps}$ of the thin layer  :
\begin{equation*}
\tilde\RR_{+}^{m;\,\eps} \times \nn= \tilde\RR_{\layer}^{m;\,\eps} \times \nn  + \mathcal{O}(\eps^{m+1}) \quad \mbox{on} \quad \Gamma^{\eps} \ ,
\end{equation*}
 see~\cite{DPP11}, and therefore  $\tilde\RR_{-}^{m;\,\eps}$ does not define an element of $\bH(\curl,\Omega)$  for smooth data. That is why the WKB expansion provides in general sharper error estimates than the multiscale expansion given in~\cite[Th. 2.9]{DPP11}.

\Bk

 \subsection{Further notes for deriving the GITC model}
 The idea for deriving the second order model with generalized impedance transmission conditions \eqref{GITC} is to write the problem satisfied by $\EE_{1}^\eps=\EE^{0}+\eps \EE^{1}$,  where $\EE^{0}$ and $\EE^{1}$ are defined  by \eqref{0E0} and \eqref{0E1}, respectively, and then to replace  $\EE_{1}^\eps$ on the left-hand sides by $\EE_{[1]}^\eps$  and to replace  $\eps\EE^{0}$ on the right-hand sides by   $\eps\EE_{[1]}^\eps$, see also for instance~\cite[Sect. 4]{Durufle2014}.

\section{Conclusion and prospects}

This work provides the first terms of a WKB expansion  for the electric field in power series of $\varepsilon$ and new second order asymptotic models with generalized impedance transmission conditions for the electric field. This work provides also mixed variational formulations for second order approximate models. A next work will be to compare numerically the accuracy of the new parameterizations $\EE^0 +\varepsilon
\EE^1$  for the electric field  as well as the second order GITC models with previous works~\cite{DPP11,Durufle2014}. In particular,  we expect that the locking phenomenon observed in~\cite{Durufle2014} with the parameterization   for the electric field in the  low frequency case does not occur anymore with the new parameterization $\EE^0 +\varepsilon \EE^1$. Finally, it will be interesting also to compare the accuracy of the new second order models with both the symmetric model  defined by~\cite[Eq. (4.8)]{Durufle2014}, 
coming from Delourme {\it  et al.}~\cite{DHJ12}, and the model defined by~\cite[Eq. (4.4)]{Durufle2014} which comes from Chun {\it  et al.}~\cite{CHH10}. 

\clearpage
\bibliographystyle{plain}
\bibliography{biblio}

\begin{thebibliography}{10}

\bibitem{ammari2014mathematical}
H.~Ammari, J.~Garnier, and L.~Giovangigli.
\newblock Mathematical modeling of fluorescence diffuse optical imaging of cell
  membrane potential changes.
\newblock {\em Quarterly of Applied Mathematics}, 72(1):137--176, 2014.

\bibitem{ammari1998effective}
H.~Ammari and S.~He.
\newblock Effective impedance boundary conditions for an inhomogeneous thin
  layer on a curved metallic surface.
\newblock {\em IEEE Transactions on Antennas and Propagation}, 46(5):710--715,
  1998.

\bibitem{Braess:2007}
D.~Braess.
\newblock {\em Finite Elements: Theory, Fast Solvers, and Applications in Solid
  Mechanics}.
\newblock Cambridge University Press, 3th edition, 2007.

\bibitem{cakoni2018nondestructive}
F.~Cakoni, I.~De~Teresa, and P.~Monk.
\newblock Nondestructive testing of delaminated interfaces between two
  materials using electromagnetic interrogation.
\newblock {\em Inverse Problems}, 34(6):065005, 2018.

\bibitem{CDFP11}
{G}. {C}aloz, {M}. {D}auge, {E}. {F}aou, and {V}. {P}{\'e}ron.
\newblock On the influence of the geometry on skin effect in electromagnetism.
\newblock {\em Comput. Methods Appl. Mech. Engrg.}, 200(9-12):1053--1068, 2011.

\bibitem{CHH10}
S.~Chun, H.~Haddar, and J.~S. Hesthaven.
\newblock High-order accurate thin layer approximations for time-domain
  electromagnetics, {P}art {II}: transmission layers.
\newblock {\em J. Comput. Appl. Math.}, 234(8):2587--2608, 2010.

\bibitem{DHJ12}
B.~Delourme, H.~Haddar, and P.~Joly.
\newblock On the well-posedness, stability and accuracy of an asymptotic model
  for thin periodic interfaces in electromagnetic scattering problems.
\newblock {\em Mathematical Models and Methods in Applied Sciences}, 13
  (23):2433--2464, 2013.

\bibitem{DPP11}
M.~Durufl{\'e}, V.~P{\'e}ron, and C.~Poignard.
\newblock Time-harmonic {M}axwell equations in biological cells---the
  differential form formalism to treat the thin layer.
\newblock {\em Confluentes Math.}, 3(2):325--357, 2011.

\bibitem{Durufle2014}
M.~Durufl{\'e}, V.~P{\'e}ron, and C.~Poignard.
\newblock Thin layer models for electromagnetism.
\newblock {\em Commun. Comput. Phys.}, 16(1):213--238, 2014.

\bibitem{EN93}
B.~Engquist and J.C. N\'ed\'elec.
\newblock Effective boundary condition for acoustic and electromagnetic
  scattering in thin layers.
\newblock {\em Technical Report of CMAP}, 278, 1993.

\bibitem{goffi2020approximate}
F.~Z. Goffi, K.~Lemrabet, and T.~Arens.
\newblock Approximate impedance for time-harmonic maxwell's equations in a non
  planar domain with contrasted multi-thin layers.
\newblock {\em Journal of Mathematical Analysis and Applications},
  489(1):124141, 2020.

\bibitem{goffi2017transfer}
F.~Z. Goffi, K.~Lemrabet, and T.~Laadj.
\newblock Transfer and approximation of the impedance for time-harmonic
  maxwell’s system in a planar domain with thin contrasted multi-layers.
\newblock {\em Asymptotic Analysis}, 101(1-2):1--15, 2017.

\bibitem{haddar2002stability}
H.~Haddar and P.~Joly.
\newblock Stability of thin layer approximation of electromagnetic waves
  scattering by linear and nonlinear coatings.
\newblock {\em Journal of computational and applied mathematics},
  143(2):201--236, 2002.

\bibitem{khelifi2016small}
A.~Khelifi and S.~Boujemaa.
\newblock Small perturbation of a surface: Full maxwell's equations.
\newblock {\em Journal of Mathematical Analysis and Applications},
  444(2):1721--1738, 2016.

\bibitem{N01}
J.-C. N{\'e}d{\'e}lec.
\newblock {\em Acoustic and electromagnetic equations}, volume 144 of {\em
  Applied Mathematical Sciences}.
\newblock Springer-Verlag, New York, 2001.
\newblock Integral representations for harmonic problems.

\bibitem{Pe09}
V.~P\'eron.
\newblock {\em Mod\'elisation math\'ematique de ph\'enom\`enes
  \'electromagn\'etiques dans des mat\'eriaux \`a fort contraste}.
\newblock PhD thesis, Universit\'e Rennes 1, France, 2009.
\newblock http://tel.archives-ouvertes.fr/tel-00421736/fr/.

\bibitem{Pe16}
V.~P{\'e}ron, K.~Schmidt, and M.~Durufl{\'e}.
\newblock Equivalent transmission conditions for the time-harmonic maxwell
  equations in 3d for a medium with a highly conductive thin sheet.
\newblock {\em SIAM Journal on Applied Mathematics}, 76(3):1031--1052, 2016.

\bibitem{stupfel2021well}
B.~Stupfel, P.~Payen, and O.~Lafitte.
\newblock A well-posed and effective high-order impedance boundary condition
  for the time-harmonic scattering problem from a multilayer coated 3-d object.
\newblock {\em Progress In Electromagnetics Research B}, 94:127--144, 2021.

\end{thebibliography}

\end{document}